\documentclass[12pt]{amsart}
\usepackage{amssymb,amsmath,amsfonts,amscd,euscript}
\usepackage{tikz}
\usepackage{tikz-cd}
\usepackage{fullpage}
\usepackage{stmaryrd}

\usepackage{hyperref}
\usepackage{todonotes}

\newcommand{\nc}{\newcommand}

\numberwithin{equation}{section}
\newtheorem{thm}{Theorem}[section]
\newtheorem*{thm*}{Theorem}

\newtheorem{prop}[thm]{Proposition}
\newtheorem{lem}[thm]{Lemma}
\newtheorem{cor}[thm]{Corollary}
\theoremstyle{remark}
\newtheorem{rem}[thm]{Remark}
\newtheorem{definition}[thm]{Definition}
\newtheorem{example}[thm]{Example}
\newtheorem{dfn}[thm]{Definition}

\nc{\gl}{\mathfrak{gl}}
\nc{\GL}{\mathsf{GL}}
\nc{\g}{\mathfrak{g}}
\nc{\gh}{\widehat\g}
\nc{\h}{\mathfrak{h}}
\nc{\la}{\lambda}
\nc{\al}{\alpha }
\nc{\be}{\beta }
\nc{\ve}{\varepsilon }
\nc{\om}{\omega }
\nc{\br}{{\bf r}}
\nc{\tr}{{\rm tr}}
\nc{\ta}{\theta}
\nc{\ch}{{\mathop {\rm ch}}}
\nc{\Tr}{{\mathop {\rm Tr}\,}}
\nc{\Id}{{\mathop {\rm Id}}}
\nc{\ad}{{\mathop {\rm ad}}}
\nc{\bra}{\langle}
\nc{\ket}{\rangle}
\nc{\pa}{\partial}
\nc{\ld}{\ldots}
\nc{\cd}{\cdots}
\nc{\hk}{\hookrightarrow}
\nc{\T}{\otimes}
\nc{\gr}{\mathrm{gr}}
\nc{\ov}{\overline}

\nc{\cO}{\mathcal O}
\nc{\msl}{\mathfrak{sl}}
\nc{\mgl}{\mathfrak{gl}}
\nc{\U}{\mathrm U}
\nc{\V}{\EuScript V}
\nc{\cL}{\mathcal{L}}

\newcommand{\bZ}{{\mathbb Z}}

\newcommand{\bP}{{\mathbb P}}

\newcommand{\fh}{{\mathfrak h}}

\newcommand{\fg}{{\mathfrak g}}

\newcommand{\fb}{{\mathfrak b}}

\newcommand{\fn}{{\mathfrak n}}
\newcommand{\A}{\EuScript{A}}

\newcommand{\ba}{{\bf a}}

\nc{\cat}{\mathcal{C}}
\nc{\RepGL}[2]{{\mathsf{Rep}}(\GL_{#1})^{(#2)}}
\nc{\RepGLt}[2]{{\mathsf{Rep}}(\mgl_{#1}[t])^{(#2)}}
\nc{\RepSL}[2]{{\mathsf{Rep}}(\SL_{#1})^{(#2)}}
\nc{\RepSLt}[2]{{\mathsf{Rep}}(\msl_{#1}[t])^{(#2)}}
\nc{\Rep}{\mathsf{Rep}}
\nc{\de}{{\text -}}

\begin{document}

\title[Semi-infinite Pl\"ucker relations and arcs over toric degeneration]
{Semi-infinite Pl\"ucker relations and arcs over toric degeneration}

\author{Ievgen Makedonskyi}
\address{Ievgen Makedonskyi:\newline
 Centre of Perspective Studies, SkolTech;
Nobel st., 1, Moscow, Russia}
\email{makedonskii\_e@mail.ru}

\begin{abstract}
We study the algebra of Weyl modules in types $A$ and $C$ using the methods of
arcs over toric degenerations and functional realization of
dual space. We compute the generators and relations
of this algebra and construct its basis.
\end{abstract}

\maketitle

\section*{Introduction}

The goal of this paper is to study the homogenous coordinate ring of the semi-infinite
flag variety (see \cite{FiMi}). This ring is graded by the cone of dominant weights and its homogenous
components are the dual global Weyl modules.

We consider a simply-connected simple algebraic group $G$ over an algebraically closed field
$\Bbbk$ of characteristic zero. We are interested in the corresponding current group
$G[[t]]$ and the corresponding current Lie algebra $\fg[t]:=\fg \otimes \Bbbk[t]$.
The group $G[[t]]$ is the group of $\Bbbk[[t]]$-points of $G$, which as variety is the arc (or in the other terms $\infty$-jets) scheme
of $G$.

The main representation theoretic objects of this paper are Weyl modules (see \cite{CP,CL1,CL2}).
Their role in the representation theory of the current algebras is similar to the role of
the irreducible representations in the representation theory of semisimple Lie algebras (see \cite{Kat,FMO,BF2,I}).
A Weyl module $\mathbb{W}_{\lambda}$ is a cyclic module with a generating highest weight vector $v_\lambda$
defined by the set of relations which essentially mean that this module is integrable and all weights of its vectors
are less or equal to $\lambda$.

We consider a unipotent subgroup $U \subset G$ and a maximal torus $H \subset G$.
The semi-infinite flag variety is the variety $G[[t]]/(U[[t]] \cdot H)$.
The variety $G[[t]]/U[[t]]$ is its affine multicone. We study the homogenous coordinate ring of
$G[[t]]/(U[[t]] \cdot H)$ which is equal to the coordinate ring of the affine variety $G[[t]]/U[[t]]$.
This ring is isomorphic to the ring of dual Weyl modules (\cite{BF1,BF2,BF3}):
\[\mathbb{W}=\bigoplus_{\lambda \in P_+} \mathbb{W}_\lambda^*.\]

The algebra $\mathbb{W}$ is generated by the sum of dual fundamental Weyl modules. Our goal is to write down the defining relations of this ring
and to give its basis. We call them the semi-infinite Pl\"ucker relations.
They were computed for type $A$ in \cite{FeMa2}. It is relatively easy to write down the set of
quadratic relations satisfied by
generators of $\mathbb{W}$. However it's hard to prove that these relations generate the whole ideal
of relations in $\mathbb{W}$. In \cite{FeMa2} this claim was proven using the fusion product structure on
local Weyl modules \cite{CL1} and this proof is very technical. Here we propose another method to prove the nonexistence of additional relations.
This method is the current version of the toric degeneration method due to Gonciulea-Lakshmibai \cite{GL}.
We study the algebra $\mathbb{W}$ inside the algebra of functions on the current group. This algebra
is usually much easier then the algebra $\mathbb{W}$ and has a structure of $\fg[t]\de \fg[t]$ bimodule.
Then the algebra $\mathbb{W}$ is the ring of invariants of $\Bbbk[G[[t]]]$ with respect to the right action of
$\fn_+[t]$. We study the current version of the toric degeneration which is the ring of lower terms
of elements $\mathbb{W}$ with respect to some order on elements of $\Bbbk[G[[t]]]$.
This algebra has obviously the graded dimension less then or equal to the graded dimension of $\mathbb{W}$.
Thus if we prove some lower bound for the graded dimension of the current version of toric degeneration then
we prove the same lower bound for $\mathbb{W}$. We use this method to prove that the algebra $\mathbb{W}$
surjects to the quotient by the quadratic semi-infinite Pl\"ucker relations.

We apply these methods  in types $A$ and $C$.
The good property of these types which we use is that the coordinate ring of $G[[t]]/U[[t]]$ is the
quotient of the coordinate ring of the arc scheme over $G/U$ by the ideal of nilpotents. This is equivalent to the
equality $W_{\omega_i} \simeq V_{\omega_i}$ for local Weyl modules of each fundamental weight $\omega_i$.
Another good property is that the coordinate rings of simply-connected
groups of these types are easy to describe. However we plan to apply our methods in a more general situation.
In particular we essentially study not the flag variety, but its toric degeneration. And
the same methods can be applied to the study of more general toric varieties.

In types $A$ and $C$ we obtain the character formulas for Weyl modules. In type $A$ this formula was proven in \cite{FeMa2}, in
type $C$ it is new.

More precisely we prove the following Theorem.

\begin{thm}
The algebra $\mathbb{W}$ in type $SP_{2n}$ is generated by allowed minors $m_I^{(k)}$ with
defining relations \eqref{SPsnakeplueckerequation} and its basis is \eqref{SPnBasisMonomials}.
In particular the character or the Weyl module $\mathbb{W}(\sum \lambda_p \omega_p)$ is the following:
\[\sum_{\br} \frac{q^{\sum_{I \prec J}k(I,J)r_Ir_J}}
{\prod_{I} (q)_{r_I}}(\prod_{I}\prod_{k=1}^{|I|} e^{\varepsilon_{i_k}})^{r_I},\]
where the summation is on $\br$ such that $\sum_{|I|=p} r_I=\lambda_p$.
\end{thm}

The structure of the paper is the following. In Section $1$ we explain the classical construction of Pl\"ucker equations and
toric degeneration of flag variety. In Section $2$ we recall the basic facts from the theory of current algebras and Weyl modules.
In Section $3$ we recall the construction of semi-infinite Pl\"ucker relations, construct an analogue of the toric degeneration
construction and give a simple proof of the fact that the  semi-infinite Pl\"ucker relations generate the ideal of
relations in $\mathbb{W}$. In Section $4$ we reduce the case of $SP_{2n}[[t]]$ to the case of $SL_{2n}[[t]]$.

\section{Finite dimensional case}
In this paper all groups, varieties and algebras are defined over algebraically closed field $\Bbbk$ of characteristic zero.

\subsection{Algebra of function on the affine cone of the flag variety}
Let $\fg$ be a finite dimensional simple Lie algebra, $G$ be the corresponding connected Lie group. We take a Borel subalgebra $\fb$.
 Let $P$ be the lattice of
weights of the group $G$, $P_+$ be the cone of dominant weights.
For $\lambda\in P_+$ we denote the corresponding simple (left) module $V_\lambda$ and let $V_\lambda^o$ be the corresponding right module.
Let $W$ be the Weyl group of the Lie algebra $\fg$ and $w_0$ be the longest element of this Weyl group.

Consider the algebra $\Bbbk[G]$ of algebraic functions on the group $G$. It has the natural structure of $G\de G$ bimodule. Therefore it has the structure of
$\fg\de \fg$ bimodule. The structure of this bimodule is described by the classical Peter-Weyl theorem.
\begin{thm}
\[\Bbbk[G]\simeq \bigoplus_{\lambda \in P_+}V_\lambda^* \otimes_{\Bbbk}V_{\lambda}^{*o}.\]
\end{thm}

Then the invariant space of this module with respect to the right action of $\fn_+$ is as follows.
\begin{cor}
\[\Bbbk[G]^{\fn_+}\simeq \bigoplus_{\lambda \in P_+}V_\lambda^*.\]
\end{cor}
This is the coordinate ring of the affine cone over flag variety.
The right action of the Cartan subalgebra $\fh$ preserves this invariant subalgebra and makes it $P_+$ graded algebra, $deg V_\lambda^*=\lambda$.

In many cases the structure of a group is easier than the structure of the corresponding
 flag variety and it is convenient to study the coordinate ring of the flag variety
inside the coordinate ring of the group.
The goal of this section is to explain some classical results on the algebra $\mathcal{M}=\bigoplus_{\lambda \in P_+}V_\lambda^*$ of functions on the affine cone of the flag variety.

We consider two dominant weights $\lambda, \mu \in P_+$. Then there is an injection:
\[V_{\lambda+\mu}\hookrightarrow V_\lambda \otimes_{\Bbbk} V_\mu.\]
Dualizing this injection we obtain the surjection:
\[V_\lambda^* \otimes_{\Bbbk} V_\mu^*\twoheadrightarrow V_{\lambda+\mu}^*.\]
This is the multiplication map in $\mathcal{M}$. In particular,
surjectivity of this map means that $\mathcal{M}$ is generated by $\bigoplus V_{\omega_i}^*$,
where the weights $\omega_i$ are fundamental.
Moreover this algebra is quadratic and even Koszul (see for example \cite{BC}).
In the following subsections we give a proof the quadraticity in some concrete
cases and describe the defining relations of this algebra.

\subsection{$SL_n$-case}
Consider the special linear group $SL_n$ of $n \times n$ matrices with determinant $1$. By definition the coordinate ring of this group is the following.
We consider $n\times n$ generators $z_{uv}$. Then:
\begin{equation}\label{SLnAlgebraofFunctions}
\Bbbk[SL_n]\simeq \Bbbk[z_{uv}]\left/\Bigg \langle det\left(
                                                 \begin{array}{ccc}
                                                   z_{11} & \cdots & z_{1n} \\
                                                   \vdots & \ddots & \vdots \\
                                                   z_{n1} & \cdots & z_{nn} \\
                                                 \end{array}
                                               \right)-1
\Bigg \rangle\right..
\end{equation}

Here the left and the right actions of $\msl_n$ are defined from the obvious matrix multiplication, i. e.
for an unit matrix $E_{ij}$ we have:
\[E_{ij} z_{uv}=\begin{cases} z_{jv}, ~\text{if}~ i=u; \\
0, ~\text{otherwise}.
\end{cases}\]

\[ z_{uv}E_{ij}=\begin{cases} z_{ui}, ~\text{if}~ j=v; \\
0, ~\text{otherwise}.
\end{cases}\]

We are interested in the precise description of the algebraic functions on the special linear group invariant under the right action of the subgroup of unitriangular matrices (or in other words invariant under the action of Lie algebra $\fn_+$).
Let $I\subset \{1, \dots, n\}$ be a nonzero proper subset, let $i_1<i_2<\dots <i_k$ be the elements of $I$.
We denote:
\[m_{I}=det\left(
                                                 \begin{array}{ccc}
                                                   z_{i_11} & \cdots & z_{i_1k} \\
                                                   \vdots & \ddots & \vdots \\
                                                   z_{i_k1} & \cdots & z_{i_kk} \\
                                                 \end{array}
                                               \right).\]

The following result is classical (see \cite{MS}):
\begin{thm}\label{MinorsAlgebraSLn}
\[\Bbbk[SL_n]^{\fn^+} =\Bbbk[m_I]_{I \subset \{1, \dots, n\}, |I|\neq 0, n}.\]
\end{thm}

\subsubsection{Construction of monomial order}
We construct the order on monomials in minors which we use throughout the paper.
\begin{dfn}\label{monomialOrder}
We define the partial monomial order on monomials in variables $m_{I}$ in the following way.
Let $I^1, \dots, I^a; J^1, \dots, J^b \subset \{1, \dots, n\}$ be
two tuples of proper nonempty subsets, $|I^1|\geq \dots \geq |I^a|$, $|J^1|\geq \dots \geq |J^b|$.
If $a>b$, then we write
\[m_I:=m_{I^1}\dots m_{I^a}\succ m_{J^1}\dots m_{J^b}=:m_J.\]
Assume that $a=b$ and $(|I^1|, \dots ,|I^a|)>(|J^1|, \dots , |J^a|)$ in the lexicographic order, i. e. $\exists c:$
$|I^c|>|J^c|$, $|I^{c'}|=|J^{c'}|$ for all $c'<c$. Then we write $m_{I^1}\dots m_{I^a}\succ m_{J^1}\dots m_{J^a}$.

Assume that $\forall c:$ $|I^c|=|J^c|$ (the case of equal Young diagrams).
We consider the n-tuple $wt (m_I)=(|\{c|1 \in I^c\}|, \dots, |\{c|n \in I^c\}|)$, we call it the weight of an element. We denote by $wt_d(m_I)$
the $d$-th component of $wt (m_I)$. Then if $\exists d:$ $wt_d(m_I)>wt_d(m_J)$, $\forall  d'>d:$ $wt_{d'}(m_I)=wt_{d'}(m_J)$, then
$m_I\succ m_J$. In words, if the number $n$ appears in $I$ more times than in $J$, then $m_I\succ m_J$, if they equal then
we compare the number of $n-1$'s etc. (This step is necessary for type $C$ purposes, for type $A$ purposes it is possible
to define an easier order).

Consider now the case of equal Young diagrams and equal weights.
We fix a number $l$, $1 \leq l \leq n-1$. Define the truncation map. For $L \subset \{1, \dots, n\}$, $|L| \geq l$, we put
\[\tr^{l-1}(L):=\{i_{l}, \dots,i_{|L|}\},\]
i. e. this operation deletes from $L$ its $l-1$ smallest elements and is $(l-1)$st power of the operator $\tr$ deleting the smallest element of the set.
For $I=(I_1, \dots, I_a)$ define
\[\tr I:=(\tr I_1, \dots, \tr I_a).\]

Assume that  $\tr^l I=\tr^l J$ for some $l$ and $\tr^{l-1} I\neq \tr^{l-1} J$. Then if the multiplicity of
$c$ in $\tr^{l-1}I$ is greater than the multiplicity of
$c$ in $\tr^{l-1}J$ and for all $c'<c$ the  multiplicity of
$c'$ in $\tr^{l-1}I$ equal to the multiplicity of
$c'$ in $\tr^{l-1}J$ we put $m_I \succ m_J$.

Assume that all these multiplicities are equal, i. e. $wt(\tr^{l-1}I)=wt(\tr^{l-1}J)$.
We order sets  $U \subset \{l+1, \dots, n\}$ in the lexicographic way, i. e. for
$U=\{u_1,\dots, u_{|U|}\}$, $V=\{v_1,\dots, v_{|V|}\}$, $u_p<u_{p+1}$, $v_p<v_{p+1}$ we write
$U<V$ if $u_1=v_1, \dots, u_{c-1}=v_{c-1}, u_c<v_c$ or $u_1=v_1, \dots, u_{|V|}=v_{|V|}, |U|>|V|$.
For example, the empty set is the largest.
Note that $\tr^{l-1}(L)=\tr^{l}(L) \cup \{a\}$ for some number $a$. We put
$p_{a,U}(I)=|\{I^c| \tr^{l-1}I^c=U \cup \{a\}, \tr^{l}I^c=U\}|$.
Then we write $m_I\succ m_J$ if $\exists a, U:$ $p_{a,U}(I)>p_{a,U}(J)$ and $\forall U'>U:$ $p_{a,U'}(I)=p_{a,U'}(J)$,
$\forall a'<a, U':$ $p_{a',U'}(I)=p_{a',U'}(J)$.

This is the linear order on monomials $m_I$. It is easy to see that this order is indeed a monomial order, i. e.
if $m_I\preceq m_J$, then $m_Im_K\preceq m_Jm_K$.
\end{dfn}

Such complicated construction of the order is needed for the purposes of current algebra. We use the same order on the tuples of sets
$(I^1, \dots, I^a)$. I. e. we define $(I^1, \dots, I^a) \succ (J^1, \dots, J^b)$ iff $m_{I^1} \dots m_{I^a} \succ m_{J^1} \dots m_{J^b}$.
On the came way we order the sets of numbers $\br:=(r_I)$, $r_I\in \mathbb{N}\cup \{0\}$, where
$r_I:=|\{p|I^p=I\}|$.

\begin{rem}
The same construction of the order can be applied for $I^1, \dots, I^a; J^1, \dots, J^b $ be subsets of arbitrary linearly ordered set.
We use the ordered set $\{1<\bar 1 < \dots < n<\bar n\}$ in symplectic case. Any monotonous map on linearly ordered sets
gives the monotonous map on monomials in minors. Later we use the monotonous bijection between
$\{1,\bar 1, \dots, n,\bar n\}$ and $\{1,2,\dots,2n\}$ to identify the monomials in minors of type $C_n$ with monomials in minors of
type $A_{2n}$.
\end{rem}

\begin{example}
\[m_{123}m_{46}m_{1}m_{1}\succ m_{78}m_{67}m_{36}m_{45},\]
because we have the inequality on lengths of sets of indices, $(3,2,1,1)>(2,2,2,2)$.
\[m_{126}m_{15}m_{1}m_{1}\succ m_{123}m_{46}m_{3}m_{1},\]
because in both monomials $6$ appears $1$ time, $5$ appears in the first monomial $1$ time and $0$ times in the second.
\[m_{1268}m_{157}m_{1}m_{1}\succ m_{2468}m_{467}m_{3}m_{1},\]
because $tr^2$ of both sides coincides, $tr(1268,157,1,1)=(268,57)>(468,67)=tr(2468,467,3,1)$ because
$2$ appears in $(268,57)$ once and doesn't appear in $(468,67)$.
\[m_{1268}m_{157}m_{1}m_{1}\succ m_{2568}m_{127}m_{3}m_{1},\]
because  $tr^2(1268,157,1,1)=(68,7)=tr^2(2568,127,3,1)$,  $tr(1268,157,1,1)=(268,57)$, $tr(2568,127,3,1)=(568,27)$,
we have $\{68\}<\{7\}$, and in the first monomial the larger number $5$ is attached to the larger set $\{7\}$, in the second
monomial this number is added to the smaller set $\{68\}$. On the other terms, in the second case "snake" on
the first two columns is longer.
\end{example}

\subsubsection{Pl\"ucker relations}
Now we want to describe the ideal of relations satisfied by the elements $m_I$ and to give a linear basis of the algebra $\mathcal{M}$.
We take two sets of numbers $I, J \subset \{1, \dots, n\}$, $I=\{i_1<i_2<\dots <i_k\}$, $J=\{j_1<j_2<\dots <j_l\}$.
We call $I\leq J$ if $k \geq l$ and $i_s \leq j_s$, $1 \leq s \leq l$. This gives a partial order on the set of nonempty proper subsets
of $\{1, \dots, n\}$.
Assume that $I, J$ are uncomparable and $|I|=k>l=|J|$ or $k=l$, $i_k=j_k,i_{k-1}=j_{k-1}, \dots, i_{r+1}=j_{r+1}, i_{r}<j_{r}$ for some $r$.
Define $s=max\{s'|i_{s'}>j_{s'}\}$.
Then the set of numbers $A=\{j_1, \dots, j_s,i_s, \dots, i_k\}$ is strictly increasing. Fix an $s$-element subset $B \subset A$.
Let $B=\{b_1<\dots <b_s\}$; $A\backslash B=\{c_1<\dots <c_{k-s}\}$. Then the chain of numbers $(b_1, \dots, b_s,c_1, \dots, c_{k-s})$ is the
permutation of the chain $(j_1, \dots, j_s,i_s, \dots, i_k)$. We denote by $sign(B)$ the sign of this permutation (i. e. $1$ for even and $-1$
for odd one).

\begin{prop}\label{plueckerFinTypeA}(See \cite{MS})
For any uncomparable pair $I, J$ there exists the following quadratic relation satisfied by minors:
\[\sum_{B\subset A, |B|=s}sign(B)m_{(I\backslash A)\cup B}m_{(J\backslash A)\cup(A \backslash B)}=0.\]
The leading monomial of the left hand side of this equality is $m_Im_J$.
\end{prop}
\begin{proof}
The function
\[\sum_{B\subset A, |B|=s}sign(B)m_{(I\backslash A)\cup B}m_{(J\backslash A)\cup(A \backslash B)}\]
is multilinear and anti-symmetric in $i+1$ vectors $(z_{a,1}, \dots, z_{a,i}), a \in A$ of length $i$. Therefore it is identically equal to zero.
By the construction we know that $A \cap J$ is the set of smallest elements of $J$. Therefore $m_I$ is the leading term.
\end{proof}
\begin{cor}
\[\Bbbk[m_I]_{I \subset \{1, \dots, n\}, |I|\neq 0, n}=\sum_{I_1\leq I_2 \leq \dots\leq I_k}\Bbbk m_{I_1}m_{I_2}\dots m_{I_k}.\]
I. e. the ring $\Bbbk[m_I]_{I \subset \{1, \dots, n\}, |I|\neq 0, n}$ is linearly generated by the products of minors with comparable indices.
\end{cor}

Next we explain the proof of the linear independence of products of minors with comparable indices. This proof uses
 the construction of toric degeneration and is due to \cite{GL}.

\begin{lem}
The subalgebra $\Bbbk[z_{ij}]_{(i,j)\neq (n,n)}$ is free.
\end{lem}
\begin{proof}
The only defining relation of the coordinate ring of the special linear group is
 \[r=det \left(\begin{array}{ccc}
                                                   z_{11} & \cdots & z_{1n} \\
                                                   \vdots & \ddots & \vdots \\
                                                   z_{n1} & \cdots & z_{nn} \\
                                                 \end{array}
                                               \right)-1.\]
 This element is linear in $m_{nn}$.
 We need to prove that $\Bbbk[z_{ij}]r\cap\Bbbk[z_{ij}]_{(i,j)\neq (n,n)}=\{0\}$.
 We put $r=r_0+r_1m_{nn}$, where $r_i \in \Bbbk[z_{ij}]_{(i,j)\neq (n,n)}$. Assume that for some $a \in \Bbbk[z_{ij}]$
 $a(r_0+r_1m_{nn}) \in \Bbbk[z_{ij}], (i,j)\neq (n,n)$. We know that $a$ is a polynomial in $z_{nn}$ with coefficients in
 $\Bbbk[z_{ij}], (i,j)\neq (n,n)$. Then $a(r_0+r_1m_{nn})$ is obviously a nonconstant polynomial in $m_{nn}$ and this completes the proof.
\end{proof}

We denote $z_{ij}<z_{i'j'}$ if $j<j'$ or $j=j'$ and $i<i'$. Then define the order $\lhd$ on monomials in $z_{ij}, (i,j)\neq (n,n)$ lexicographically in the
order of variables. For $f \in \Bbbk[z_{ij}]_{(i,j) \neq (n,n)}$ define by $lt(f)$ the leading term on $f$ in $\lhd$. We apply the operation
$lt$ to the algebra $\Bbbk[m_I]$. The algebra $lt(\Bbbk[m_I])$ is called the toric degeneration of the algebra
$\Bbbk[m_I]_{I \subset \{1, \dots, n\}, |I|\neq 0, n}$.

We denote $I_k=\{i_{k1}<\dots{i_{k|I_k|}}\}$.
The next proposition is obvious.
\begin{prop}\label{SLleading}
\[lt(m_{I_1}\dots m_{I_s})=\prod_{i_{kj} \in I_k} z_{i_{kj}j}.\]
The leading terms of monomials in $m_I$ are the monomials
\[z_{i_{11} 1}z_{i_{21} 1}\dots z_{i_{l_11}1} z_{i_{22} 1}\dots z_{i_{l_22}2}\dots\]
for $l_1\geq l_2 \geq \dots$ and $i_{kj}<i_{k,j+1}$.
They are different for different products of comparable monomials.
\end{prop}

Different leading terms are clearly linearly independent.
Therefore the products of comparable monomials are linearly independent and they form a basis
of $\Bbbk[m_I]$.

\subsection{Symplectic case.} Consider now the symplectic group $SP_{2n}$. It consists of matrices $S$ such that:
\[S \left(\begin{array}{cc}
                                                   0 & E \\
                                                   -E & 0 \\
                                                 \end{array}
                                               \right) S^t=
                                               \left(\begin{array}{cc}
                                                   0 & E \\
                                                   -E & 0 \\
                                                 \end{array}
                                               \right)\]
Enumerate the rows and the columns of $2n\times 2n$ matrices by elements $1, 2, \dots, n, \bar {1}, \dots, \bar{n}$. Then we denote
the coordinate functions on the space of $2n \times 2n$ matrices by $z_{ij}$, $i, j \in \{1, 2, \dots, n, \bar {1}, \dots, \bar{n}\}$.

The symplectic group is an affine variety with the generators $z_{ij}$ and the following defining relations:
\[\sum_{k=1}^n z_{ki}z_{\bar{k}j}-z_{\bar{k}i}z_{kj}=
\begin{cases}
  1, & \mbox{if } j=\bar{i} \\
  -1, & \mbox{if } i=\bar{j} \\
  0, & \mbox{otherwise}.
\end{cases}\]

We order the elements $\{1, 2, \dots, n, \bar {1}, \dots, \bar{n}\}$ in the following way:
\[1<\bar{1}<2<\bar{2} \dots <n<\bar{n}.\]
For a set $I =\{i_1<i_2<\dots\}\subset \{1, 2, \dots, n, \bar {1}, \dots, \bar{n}\}$, $|I| \leq n$ let
$m_I$ be the minor in the columns $1, \dots, |I|$ and
rows $i_{1}, i_2, \dots$. I. e.:
\begin{equation}\label{SPdefiningRelations}
m_{I}=det\left(
                                                 \begin{array}{ccc}
                                                   z_{i_11} & \cdots & z_{i_1k} \\
                                                   \vdots & \ddots & \vdots \\
                                                   z_{i_k1} & \cdots & z_{i_kk} \\
                                                 \end{array}
                                               \right).
\end{equation}

The following results are well known (see, for example, \cite{DC}).
\begin{prop}\label{SPInvariantsGenerators}
\[\Bbbk[SP_{2n}]^{\fn^+}=\Bbbk[m_I]_{ I \subset \{1, 2, \dots, n, \bar {1}, \dots, \bar{n}\}, |I| \leq n}.\]
\end{prop}

The minors $m_I$ obviously satisfy the relations from Proposition \ref{plueckerFinTypeA} for any uncomparable sets $I, J$.
However these minors are linearly dependent. More precisely take $I \subset \{1, 2, \dots, n, \bar {1}, \dots, \bar{n}\}$, $|I| \leq n-2$.
For any $l \in \{1, \dots, n\}$ take a multiset ${I \cup \{l, \bar{l}\}}$ and define the minor
$m_{I \cup \{l, \bar{l}\}}$. It is equal to zero if this multiset has multiple elements, because then this minor has
equal columns.
\begin{prop}\label{SPlinearDependenceMinors}
For any $I \subset \{1, 2, \dots, n, \bar {1}, \dots, \bar{n}\}$, $|I| \leq n-2$
we have:
\begin{equation}\label{eqSPlinearDependenceMinors}
m_{I \cup \{1, \bar{1}\}}+m_{I \cup \{2, \bar{2}\}}+\dots +m_{I \cup \{n, \bar{n}\}}=0,
\end{equation}
where union means the union of multisets.
\end{prop}
\begin{proof}
For $I \subset \{1, 2, \dots, n, \bar {1}, \dots, \bar{n}\}, K\subset \{1, 2, \dots, n\}, |I|=|K|$ we define:
\[m_{IK}=det\left(
                                                 \begin{array}{ccc}
                                                   z_{i_1k_1} & \cdots & z_{i_1k_{|I|}} \\
                                                   \vdots & \ddots & \vdots \\
                                                   z_{i_{|I|}k_1} & \cdots & z_{i_{|I|}k_{|I|}} \\
                                                 \end{array}
                                               \right).\]
Decomposing minors $m_{I \cup \{l, \bar{l}\}}$ by two columns (indexed by $l$ and $\bar l$) we obtain that the left hand side of
\eqref{eqSPlinearDependenceMinors} is equal to:
\[\sum_{K\subset \{1, 2, \dots, |I|+2\}, |I|=|K|, \{i,j\}=\{1, 2, \dots, |I|+2\}\backslash K} m_{IK}\sum_{l=1}^n z_{li}z_{\bar{l}j}-z_{\bar{l}i}z_{lj}\]
which is equal to $0$.
\end{proof}
Note that some of minors in \eqref{eqSPlinearDependenceMinors} are equal to $0$. They are the minors $m_{I \cup \{l, \bar{l}\}}$
if $l$ or $\bar{l}$ belongs to $I$.

\begin{definition}\label{forbidderMinors}
We call {\it forbidden} the subsets $J\subset \{1, \dots,n\}, 0< |I| \leq n$ and the minors $m_J$ such that for some $a<b$ $j_b=\bar{a}$.
The remaining subsets $J\subset \{1, \dots,n\}, 0< |I| \leq n$ and minors $m_J$ we call {\it allowed}.
\end{definition}
\begin{rem}
For each forbidden minor $m_J$ there exists a number $b$ such that $j_b=b, j_{b+1}=\bar b$.
\end{rem}
\begin{rem}
It is an easy combinatorial check that for $l \geq 2$ there is $\binom{2n}{l}-\binom{2n}{l-2}$ allowed minors $m_I$, $|I|=l$.
\end{rem}

The next Lemma essentially belongs to C. De Concini \cite{DC}. We use slightly different monomial order so
we give it with the complete proof.
\begin{lem}\label{SPForbiddenMinors}
Each forbidden minor is equivalent modulo relations \eqref{eqSPlinearDependenceMinors} to the linear combination of
allowed minors which are greater then or equal to this minor with respect to the order $\succ$.
\end{lem}
\begin{proof}
For a forbidden minor $m_J$ let $b$ be the smallest number such that $j_b=b, j_{b+1}=\bar b$.
Let $A=\{a_1, \dots, a_s\}$ be the set of numbers $a_p<b$ such that both $a_p, \bar a_p \in J$,
$C=\{c_1, \dots, c_s\}$ be the set of numbers $c_p<b$ such that both $c_p, \bar c_p \not \in J$
(clearly the cardinalities of these sets are equal).

We put $D:=A \cup C \cup \{b\}$. For each $s$-element subset $E \subset D$,
$E=\{e_1, \dots, e_s\}$ we define $\bar E:=\{\bar e_1, \dots, \bar e_s\}$.
Take the $(|J|-2)$-element subset $I:=J \backslash (A \cup \bar A\cup \{b,\bar b\})\cup (E \cup \bar E)$.
Then we consider the equality \eqref{eqSPlinearDependenceMinors} for such a subset $I$.
We have $\binom{2s+1}{s}$ equalities of such a type. Note that all these equalities
are indeed the linear combination of the following variables:
\begin{itemize}
  \item for $s+1$ element subset $F \subset D$ the minors $m_{J \backslash (A \cup \bar A\cup \{b, \bar b\})\cup (F \cup \bar F)}$, there are $\binom{2s+1}{s+1}$
  such variables;
  \item variables $m_K$ such that for some $k$ the elements greater than $k$ belong or don't belong to $K$ and $J$ together and
  $k \in K$, $k \not \in J$; for such variables we have $m_K \succ m_J$.
\end{itemize}

We consider this system of $\binom{2s+1}{s}$ equalities as the system of linear equations in $\binom{2s+1}{s+1}$ variables of the first type.
We prove that this system is nondegenerate, and hence, the monomials of the first type lie in the linear span of the monomials of the second type.

The matrix of this system is the following. Its columns are enumerated by the $s+1$ element subsets of $2s+1$ element set, and its rows
are enumerated by $s$ element subsets of the same set The element of the column $A$ and the row $B$ is equal to $1$ if $B \subset A$,
$0$, otherwise. We note that this matrix is $\mathfrak{S}_{2s+1}$ equivariant, so its kernel has to be $\mathfrak{S}_{2s+1}$
invariant subset. It is easy to check that the $\binom{2s+1}{s}$ dimensional module spanned by $s+1$ element subsets of $2s+1$ element set
splits to $s+1$ irreducible representations and all these representations don't lie in the kernel of the matrix.
\end{proof}
Thus the ring $\Bbbk[SP_{2n}]^{\fn^+}$ is generated by the allowed minors and the product of any pair
of uncomparable allowed minors can be expressed as the linear combination of the products of
pairs of the comparable minors (because they satisfy the relations from Lemma \ref{plueckerFinTypeA}).

Our next goal is to give a proof of the linear independence of products of comparable allowed minors.
We will use the method of  toric degeneration analogous to the case of $GL_n$.

\begin{lem}\label{idealRelationsIntersection}
The defining ideal of $SP_{2n}$ has a trivial intersection with the polynomial subalgebra $\Bbbk[z_{uv}]_{v \in \{1, 2, \dots, n\},
(u,v) \neq (\bar k, l), k < l}$.
\end{lem}
\begin{proof}
We denote by $\mathcal{F}$ the field $\Bbbk(z_{uv})_{v \in \{1, 2, \dots, n\},(u,v) \neq (\bar k, l), k < l}$.
Consider the $n(n-1)/2$ relations
\begin{equation}\label{SPnrowsEquation}
\sum_{k=1}^n z_{ku}z_{\bar{k}v}-z_{\bar{k}u}z_{kv}=0
\end{equation}
for $u,v\in\{1,2,\dots, n\}$. These relations are the subset of relations \eqref{SPdefiningRelations}
containing only $z_{uv}, v \in \{1, \dots, n\}$.
They are linear in variables $z_{\bar{k}v}$. Consider the variables $z_{\bar{k}l}$, $k<l$. It is easy to see that
the system \eqref{SPnrowsEquation} is triangular with respect to these $n(n-1)/2$ variables.
Therefore these relations as a system of equations on $z_{\bar{k}l}$, $k<l$, can be rewritten in the form
\[z_{\bar{k}l}-a_{\bar{k}l},~a_{\bar{k}l} \in \mathcal{F}.\]
Thus they generate a proper ideal in the polynomial ring $\mathcal{F}[z_{\bar{k}l}]=\mathcal{F}[z_{\bar{k}l}-a_{\bar{k}l}]$.
Hence this ideal has trivial intersection with the field $\mathcal{F}$ and thus with the ring
$\Bbbk[z_{uv}]_{v \in \{1, 2, \dots, n\},(u,v) \neq (\bar k, l), k < l}$ as well.
Moreover we have:
\[\mathcal{F}[z_{\bar{k}l}]/\langle z_{\bar{k}l}-a_{\bar{k}l} \rangle\simeq \mathcal{F}.\]
Therefore the elements $z_{uv}$, $(u,v)=(\bar k,l), k<l$ belong to $\mathcal{F}$.

Take now an algebra $\mathcal{F}[z_{u \bar v}]_{u \in \{1, \bar 1, \dots, n ,\bar n\}, v \in \{1, \dots, n\}}$.
Consider the remaining relations \eqref{SPdefiningRelations}.
They form an ideal in this ring. We need to prove that this ideal is proper, i. e. it contains a point.
However the vectors $(z_{uv})$, with fixed $v \in \{1, \dots, n\}$ span a Lagrangian subspace $L$
in $2n$ dimensional space over $\mathcal{F}$. Thus there is a Lagrangian subspace $L$ spanned by
$(z_{uv})$, with fixed $v \in \{\bar 1, \dots, \bar n\}$, such that the vectors $(z_{uv})$, $v \in \{\bar 1, \dots, \bar n\}$ are symplectic dual to
$(z_{uv})$, $v \in \{1, \dots, n\}$. This gives a needed point. This completes the proof.
\end{proof}

\begin{thm}\label{SPleading}
Leading monomials for different products of comparable allowed minors are different.
\end{thm}
\begin{proof}
The leading monomial of an allowed minor $m_I$ is equal to
$z_{i_1 1} \dots z_{i_{|I|},|I|} \subset \Bbbk[z_{uv}]_{v \in \{1, 2, \dots, n\},(u,v) \neq (\bar k, l), k < l}$.
We identify the leading monomials of type $SP_{2n}$ with some leading monomials of type $SL_{2n}$.
Renaming the indices $(u,v) \mapsto (2u-1,v)$,
$(\bar u,v) \mapsto (2u,v)$ for $u,v \in \{1, 2, \dots, n\}$, the comparable minors come to comparable.
Thus this theorem is the partial case of
Theorem \ref{SLleading}.
\end{proof}
\begin{cor}
Different products of comparable allowed minors are different. In particular the products of comparable allowed minors form a basis
of $\Bbbk[G/U]$. They are called symplectic Young tableaux in \cite{DC}.
\end{cor}

\section{Current algebras and Weyl modules}
In this section we recall some properties of current groups and algebras.
We call the current Lie algebra the Lie algebra $\fg \otimes \Bbbk[t]$
which is a maximal parabolic subalgebra of the untwisted affine Kac-Moody Lie algebra attached to $\fg$.
Consider a finite dimensional Lie group $G$. We call the current group $G[[t]]$ the group of $\Bbbk[[t]]$-points of
the group $G$. If $\fg$ is the Lie algebra of the Lie group $G$ then $\fg \otimes \Bbbk[t]$
is the Lie algebra of the group $G[[t]]$.

The variety $G[[t]]$ is in the other terms called the variety of arcs over $G$ (\cite{Na}). For an affine variety the arc variety
 has the following description.

Assume that an affine variety $V$ has the coordinate ring
\[\Bbbk[x_1\dots, x_r]/\langle p_1(x_1\dots, x_r), \dots, p_l(x_1\dots, x_r) \rangle.\]
Then the arc variety has the following coordinate ring.
We take the set of variables $x_i^{(k)}$, $i=1, \dots, r; k =0, 1,\dots$ and define the formal series:
\begin{equation}\label{powerseriesJetVariables}
x_i(s)=\sum_{k=0}^\infty x_i^{(k)}.
\end{equation}

Then the coordinate ring of the arc variety is
\[\Bbbk[x_i^{(k)}]/\langle p_1(x_1(s),\dots, x_r(s)), \dots, p_l(x_1(s),\dots, x_r(s)) \rangle,\]
where $p_1(x_1(s)\dots, x_r(s))$ means the set of all $s$-coefficients of this polynomial. (See for details \cite{FBZ})

\begin{prop}\label{JetGroup}
The arc scheme of an algebraic group is reduced, i. e. its coordinate ring has no nilpotents.
\end{prop}
\begin{proof}
This scheme is the group scheme (of $\Bbbk[[t]]$ points). Therefore it is reduced (\cite{Oo}).
\end{proof}

We are interested in the following representations of the current Lie algebra.

\begin{definition}\label{weylmodules}\cite{CP}
Let $\lambda\in P_+$. Then the global Weyl module
$\mathbb{W}(\lambda)$ is the cyclic $\fg \otimes \Bbbk[t]$ module with a generator $v_\lambda$ and the following
defining relations:
\begin{gather}
\label{weylvanishing1}
(e_{\alpha}\otimes t^k) v_{\lambda}=0, \ \alpha \in \Delta_+, ~k \geq 0;\\
\label{weylbound1}
(f_{-\alpha}\otimes 1)^{\langle \alpha^\vee, \lambda \rangle+1} v_{\lambda}=0, \ \alpha \in \Delta_+.
\end{gather}
Local Weyl modules ${W}(\lambda)$ are defined by previous conditions and one additional condition:
\begin{equation}
\label{weylvanishing0}
h\T t^k v_{\lambda}=0 \text{ for all } h\in\fh, k>0.\\
\end{equation}
\end{definition}

Weyl modules are the natural analogues of finite-dimensional simple $\fg$-modules $V(\lambda)$. They are graded by the degree of $t$:
\[\mathbb{W}(\mu)=\bigoplus_{k=0}^\infty \mathbb{W}(\mu)^{(k)}\]
with finite-dimensional homogeneous components. Therefore we can define the restricted dual module:
\[\mathbb{W}(\mu)^*=\bigoplus_{k=0}^\infty (\mathbb{W}(\mu)^{(k)})^*.\]

Dual global Weyl modules have the following properties.

\begin{lem}\label{cocyclic}
$\mathbb{W}(\mu)^*$ is cocyclic, i. e. there exists an element (cogenerator) $v^*\in\mathbb{W}(\mu)^*$
such that for any element $u\in\mathbb{W}(\mu)^*$
there exists an element $f \in U(\fg \otimes \mathbb{K}[t])$ such that $fu=v^*$. The set of cogenerators coincides
with $({\mathbb{W}(\mu)^{(0)}})^*\simeq V(\mu)^*$ .
\end{lem}
\begin{proof}
This is a direct consequence of the fact that Weyl module is cyclic.
\end{proof}

For a dominant weight $\lambda=\sum_{k=1}^{{\rm rk}(\fg)} m_k\omega_k$ we define
\[(q)_\lambda=\prod_{k=1}^{{\rm rk} \fg}\prod_{i=1}^{m_k}(1-q^i).\]

Each Weyl module is graded by the $\fh$-weights and by $t$-degree. For any $\fg\otimes \mathbb{K}[t]$-module $U$ with such a grading
let $U(\nu,m)$, $\nu\in\fh^*$, $m\in\bZ$ be the weight space of the corresponding weight.
\begin{definition}
\[ \ch U=\sum_{\nu,m}\dim U(\nu,m)x^{\nu}q^m.\]
\end{definition}

\begin{prop}\cite{CFK,N}\label{localglobal}
\[\ch \mathbb{W}(\mu)=\frac{\ch {W}(\mu)}{(q)_\mu}.\]

Moreover, the algebra of endomorphisms of $\mathbb{W}(\mu)$ is a polynomial algebra generated by
$h_{\alpha_i}t^r, 1 \leq r \leq \langle \mu, \alpha_i \rangle$. We denote this algebra by $\A_{\mu}$.
This makes $\mathbb{W}(\mu)$ the $\fg\de\A_{\mu}$ bimodule.
\end{prop}

\begin{lem}\label{weylinjection}\cite{Kat}
The $\fg \otimes \mathbb{K}[t]$-submodule of $\mathbb{W}(\la)\otimes\mathbb{W}(\mu)$
generated by $v_{\lambda}\otimes v_{\mu}$ is
isomorphic to $\mathbb{W}( \lambda+\mu)$.
\end{lem}

\begin{cor}\label{weylmodulering}
There exists a surjection of dual Weyl modules:
\begin{equation}\label{dualweylproduct}
\mathbb{W}^*(\la)\otimes\mathbb{W}^*(\mu)\twoheadrightarrow\mathbb{W}^*(\la+ \mu),
\end{equation}
inducing the structure of associative and commutative algebra on the space $\bigoplus_{\lambda \in P_+}\mathbb{W}^*(\lambda)$.
We denote this algebra by $\mathbb{W}=\mathbb{W}(\fg)$.
\end{cor}

\begin{rem}
The algebra $\mathbb{W}$ is an analogue of the algebra $\bigoplus_{\lambda \in X_+}{V}(\lambda)^*$.
\end{rem}

The following proposition is a direct consequence of Corollary \ref{weylmodulering}.

\begin{prop}\label{fundamentalGenerated}
$\mathbb{W}$ is generated by the space $\bigoplus_{k=1}^{{\rm rk}(\mathfrak{g})}\mathbb{W}(\omega_k)^*$.
\end{prop}

The structure of this algebra is complicated and we want to study this algebra using the algebra $\Bbbk[G[[t]]]$ of
functions on the current algebraic group.

Analogously to the usual (left) Weyl modules on can define the right Weyl modules $\mathbb{W}(\lambda)^{o}$ and
dual right Weyl modules $\mathbb{W}(\lambda)^{o*}$. These  modules are $\A_\lambda\de \fg$ bimodules.
Thus we can define the $\fg\de\fg$ bimodules
\[\mathbb{W}(\lambda)^{*}\otimes_{\A_\lambda}\mathbb{W}(\lambda)^{o*}.\]

We have the following Peter-Weyl theorem for current algebras.
\begin{thm} \cite{FKM}
For a simply-connected group $G$ there is a filtration on $\Bbbk[G[[t]]]$ with subquotients $\mathbb{W}(\lambda)^{*}\otimes_{\A_\lambda}\mathbb{W}(\lambda)^{o*}$.
The components of this filtration $F_\lambda$ are cogenerated as $\fb\de \fb$ bimodules
by elements $v_{-w_0(\lambda)} \otimes v_{-w_0(\lambda)}$
of left and right weight $-w_0(\lambda)$ and $t$-degree $0$.
\end{thm}

Recall that $w_0$ is the longest element of the Weyl group.

\begin{cor}\label{invariantWeylModules}
\[\Bbbk[G[[t]]]^{\fn_+[t]} \simeq \mathbb{W}.\]
In particular it is generated by the components of fundamental weights $\omega_i$.
\end{cor}
\begin{proof}
Consider the space of $\fn_+[t]$-invariant vectors of right $\fh$-weights $\lambda$. For an element $u$ of this
space we have $u\fn_+\neq 0\Leftrightarrow \exists \gamma\in \Delta_+, k \geq 0, wt(u e_\gamma t^k) > \lambda$ and thus $u$ is not cogenerated by element
of weight $\lambda$.

Therefore the space of elements cogenerated by $v_{-w_0(\lambda)} \otimes v_{-w_0(\lambda)}$ and $\fn_+$ invariant is equal to
\[\mathbb{W}_\lambda^{*}\otimes_{\A_\lambda}\A_\lambda.\]
Thus it is isomorphic to $\mathbb{W}_\lambda^{*}$ as the left module. We have that
\[v_{-w_0(\lambda)} \otimes v_{-w_0(\lambda)}\times v_{-w_0(\mu)} \otimes v_{-w_0(\mu)}\mapsto v_{-w_0(\lambda+\mu)} \otimes v_{-w_0(\lambda+\mu)}.\]
Thus by cocyclicity we have that the algebra structure of $\Bbbk[G[[t]]]^{\fn_+[t]}$ coincides with $\mathbb{W}$.
\end{proof}

\section{Semi-infinite Pl\"ucker relations in type $A$}
The goal of this section is to give an easier and conceptually different proof of the results of \cite{FeMa2}.
We recall the construction of quadratic semi-infinite Pl\"ucker relations. Then we prove that they form a basis
of ideal of relations using arc varieties over toric degenerations.

\subsection{Structure of the algebra of functions on current group}

Recall \eqref{SLnAlgebraofFunctions}:
\begin{equation}
\Bbbk[SL_n]\simeq \Bbbk[z_{ij}]\left/\Biggl \langle det\left(
                                                 \begin{array}{ccc}
                                                   z_{11} & \cdots & z_{1n} \\
                                                   \vdots & \ddots & \vdots \\
                                                   z_{n1} & \cdots & z_{nn} \\
                                                 \end{array}
                                               \right)-1
 \Biggr \rangle\right..
\end{equation}

We are interested in the current group of this group, or in the other terms we consider the arc variety with the following algebra of functions.
Consider the family of variables $z_{ij}^{(k)}$, $1 \leq i,j \leq n$, $k =0, 1, \dots$.
As before we define:
\[z_{ij}(s):=\sum_{k=0}^{\infty}z_{ij}^{(k)}s^k.\]
Then we have:
\begin{equation}\label{SLnCurrentALgebraFunctions}
\Bbbk[SL_n[t]]\simeq \Bbbk[z_{ij}^{(k)}]\left/\Biggl \langle det\left(
                                                 \begin{array}{ccc}
                                                   z_{11}(s) & \cdots & z_{1n}(s) \\
                                                   \vdots & \ddots & \vdots \\
                                                   z_{n1}(s) & \cdots & z_{nn}(s) \\
                                                 \end{array}
                                               \right)-1
\Biggr \rangle\right..
\end{equation}
Here as before $\langle f(s) \rangle$ means the ideal generated by all $s$-coefficients of $f(s)$.
The right action of the current Lie algebra $\fg[t]$ on this algebra is the following.

Let $e_{ab}$ be the matrix unit, then the element $e_{ab}t^c$ acts in the following way:

\[z_{uv}^{(k)}e_{ab}t^c=\begin{cases}z_{ub}^{(k-c)}, ~\text{if}~ a=v; k-c \geq 0;\\
0, ~\text{otherwise}.
\end{cases}\]

The left action can be described in a similar way.

We study the ring of $\fn_+[t]$ invariant functions with respect to the right action. The next Proposition tells that
this ring is generated by coefficients of minors on first columns.

Let $I\subset \{1, \dots, n\}$ be a nonzero proper subset, let $i_1<i_2<\dots <i_k$ be the elements of $I$.
\begin{definition}\label{SLnCurrentMinors}
\[m_{I}(s)=det\left(
                                                 \begin{array}{ccc}
                                                   z_{i_11}(s) & \cdots & z_{i_1k}(s) \\
                                                   \vdots & \ddots & \vdots \\
                                                   z_{i_k1}(s) & \cdots & z_{i_kk}(s) \\
                                                 \end{array}
                                               \right)\]

  \end{definition}

Expanding $m_{I}(s)$ we denote its coefficients:
\[m_{I}(s)=:\sum_{k=0}^\infty m_I^{(k)}s^k.\]
\begin{prop}\label{SLnCurrentInvariants}
\[\Bbbk[SL_n[[t]]]^{\fn_+[t]}=
\Bbbk[m_I^{(k)}]_{I \subset \{1, \dots, n\}, |I|\neq 0, n,k \geq 0}\]
\end{prop}
\begin{proof}
Clearly each $m_I^{(k)}$ is invariant under $\fn_+[t]$. The Weyl module
$W_{\omega_p}$ has the graded dimension equal to the graded dimension of the linear span of these elements with $|I|=p$.
Therefore the space of invariant functions of the right weight $\omega_p$
is equal to the span of these elements. Therefore Proposition \ref{fundamentalGenerated}
and Corollary \ref{invariantWeylModules} imply the desired equality.
\end{proof}

\subsection{Semi-infinite Pl\"ucker relations}
From this moment till the end of this Section we study the algebra generated by $m_I^{(k)}$. We have proved
 that this algebra is isomorphic to the algebra of dual Weyl modules. We denote the algebra generated by $m_I^{(k)}$ by $\mathbb{W}$ too.
Denote by $\mathbb{W}[[s]]$ the algebra of formal series in variable $s$ over algebra $\mathbb{W}$.
 First we recall
 semi-infinite Pl\"ucker relations.

Recall the order on proper nonzero subsets of $\{1, \dots, n\}$ \ref{monomialOrder}.
\begin{dfn}\label{kdefinition}
Assume that for some $I\prec J$ we have the following set of inequalities:
\begin{gather*}
i_{{|J|}}\leq j_{|J|}, \dots, i_{k_1+1}\leq i_{k_1+1},\\
i_{k_1}>j_{k_1}, i_{k_1-1}\geq j_{k_1-1}, \dots, i_{k_2+1}\geq j_{k_2+1},\\
i_{k_2}<j_{k_2}, i_{k_2-1}\leq j_{k_2-1},\dots.
\end{gather*}
We define strictly decreasing sequence of elements
\[
S(I,J)=(i_{|I|}, i_{|I|-1}, \dots, i_{|J|+1},
\dots, i_{k_1+1}, i_{k_1}, j_{k_1}, j_{k_1-1}, \dots, j_{k_2}, i_{k_2}, \dots).
\]
We set
\[
k(I,J)=|S(I,J)|-|I|.
\]
\end{dfn}

We call the set $S(I,J)$ the snake. The number $k(I,J)$ measures how much uncomparable are $I$ and $J$,
$k(I,J)=0$ iff $I$ and $J$ are comparable.

These monomials satisfy the following relations (see \cite{FeMa2}).

\begin{prop}\label{propsnakeplueckerequation}
$a)$.\ For any $k' \leq k(I,J)-1$ we have the following equality in $\mathbb{W}[[s]]$:
\begin{equation}\label{snakeplueckerequation}
\sum_{A\subset S(I,J), |A|=|I \cap S|}(-1)^{{\rm sign}(A)}
\frac{\partial^{k'} m_{I \backslash S \cup A}(s)}{\partial s^{k'}} m_{I \backslash S \cup (S \backslash A)}(s)=0.
\end{equation}
$b).$\ The series
\[\frac{\partial^{k'} m_{I}(s)}{\partial s^{k'}}m_{J}(s)\]
is the leading part in the left hand side of the previous equation.
\end{prop}

Recall the monomial order from Definition \ref{monomialOrder}.
We define the order on monomials $m_{I^1}^{(k_1)}\dots m_{I^a}^{(k_a)}$ such that
\[m_{I^1}^{(k_1)}\dots m_{I^a}^{(k_a)}\succcurlyeq m_{J^1}^{(l_1)}\dots m_{J^b}^{(l_b)} \Leftrightarrow
m_{I^1}\dots m_{I^a}\succcurlyeq m_{J^1}\dots m_{J^b}.\]
In particular this order doesn't fill the upper indices and elements which differ only by the upper indices are equivalent.

We now consider the degeneration of equations $\eqref{snakeplueckerequation}$ with respect to the order $\prec$.

\begin{dfn}\label{SLnDegeneratedAlgebra}
We denote by $\widetilde{\mathbb{W}}$ the
free polynomial algebra with generators $\bar m_I^{(k)}$ modulo the relations
\begin{equation}\label{SLnDegeneratedRelations}
\frac{\partial^{k'} \widetilde m_{I}(s)}{\partial s^{k'}}\widetilde m_{J}(s)=0,~ 0 \leq k' < k(I,J).
\end{equation}

\end{dfn}

Clearly the character of the degenerate algebra is greater then or equal to the character of the algebra $\mathbb{W}$.

We use the notation $\alpha_i=\varepsilon_i-\varepsilon_{i+1}$; $\omega_k=\sum_{p=1}^k\varepsilon_p$.
In these notation the grading of an element $m_I^{(k)}$ is equal to $\sum_{p=1}^{|I|} \varepsilon_{i_p}q^k$.
Note that the semi-infinite Pl\"ucker relations are homogenous with respect to this grading.

Degenerate Pl\"ucker relations are homogenous with respect to a stronger grading. More precisely
consider the free commutative semigroup generated by nonempty proper subsets $I \subset \{1, \dots, n\}$ and the variable $\delta$.
Attach the grading $I + k \delta$ to the variable $m_I^{(k)}$. Then it is clear that the relations \eqref{SLnDegeneratedRelations}
are homogenous with respect to this grading. We use the notation $x^\delta=q$.
For a $\mathbb{N} \delta$-graded space $V=\bigoplus_{p}V_p$ we define its $q$-dimension by $\sum_p \dim{V_p}q^p$.

 Our next goal is to compute the $q$-dimension of $\mathbb{N}[I]$-homogenous components of the algebra
$\widetilde{\mathbb{W}}$ and to describe its basis.
Denote these components by $\widetilde{\mathbb{W}}_\br$ for $\br=(r_I)$.

The next Proposition is proven in \cite{FeMa2}. We give it with the complete proof because we use the methods of this proof
in the study of arcs over toric degeneration below.
Recall the notation $(q)_r=\prod_{i=1}^r (1-q^i)$.
\begin{prop}\label{SLnDegeneratedCharacter}
\begin{equation}
\ch \widetilde{\mathbb{W}}_\br=
\frac{q^{\sum_{I \prec J}k(I,J)r_Ir_J}}
{\prod_{I} (q)_{r_I}}.
\end{equation}
\end{prop}
\begin{proof}
We consider a functional realization of the dual space of
$\widetilde{\mathbb{W}}_\br$. Namely, given a linear function $\xi$ on the space
$\widetilde{\mathbb{W}}_\br$,
we attach to it the polynomial $f_\xi$ in variables $Y_{I,p}$, $I\subset \{1,\dots,n\}$, $1\leq p\le r_I$ defined as follows:
\begin{equation}\label{dual}
f_\xi=\xi \bigl(\prod_{I} r_I(Y_{I,1})\dots r_I(Y_{I,r_I})\bigr).
\end{equation}
We claim that formula \eqref{dual} defines an isomorphism between the space of functionals on
$\widetilde{\mathbb{W}}(\sum_I r_I I)$
and the space ${\rm Pol}(\sum_I r_I I)$ of polynomials $f$ in variables $Y_{I,p}$ subject to the following conditions:
\begin{itemize}
\item $f$ is symmetric in variables $Y_{I,j}$ for each $I$,
\item $f$ is divisible by $(Y_{I,j_1}-Y_{J,j_2})^{k(I,J)}$ for all $I,J,j_1,j_2$.
\end{itemize}
The first condition holds because of commutativity of multiplication in $\widetilde{W}$.
The second one comes from the relations \eqref{SLnDegeneratedRelations}. Indeed these relations mean that this space contains
polynomials in $Y_{I,p}$ which go to $0$ under the substitution $Y_{I,p}=Y_{J,q}$ with their $Y_{I,p}$-derivatives from $0$th
to $k'(I,J)-1$st.

We note that the $q$-dimension on $\widetilde{\mathbb{W}}_\br$
 is now translated into the counting of the total degree in all variables $Y_{I,p}$.
Now the $q$-character of the space  ${\rm Pol}(\sum_I r_I I)$ is given by
\begin{equation}\label{Polch}
\frac{q^{\sum_{I\prec J}k(I,J)r_{I}r_{J}}}
{\prod_{\sigma} (q)_{r_\sigma}}.
\end{equation}
Indeed, $(q)_r^{-1}$ is the character of the space of symmetric polynomials in $r$ variables and the factor
$\prod_{I\prec J}(Y_{I,p_1}-Y_{\tau,p_2})^{k(I,J)}$ produces the numerator of the above formula.
\end{proof}

\begin{cor}\label{corssst}
The homogeneous component $\widetilde{\mathbb{W}}_\br \subset \widetilde{\mathbb{W}}$ has the basis consisting of monomials of the form
\begin{equation}\label{SLnBasisMonomials}
\prod_{I\subset\{1,\dots,n\}} m_I^{(l_{1,I})}\dots m_I^{(l_{r_I,I})},\quad 0\le l_{1,I}\le \dots\le l_{r_I,I}
\end{equation}
such that $l_{1,I}\geq \sum_{J \prec I} k(I,J) r_J$.
\end{cor}
\begin{proof}
We note that the character of the set of monomials \eqref{SLnBasisMonomials} is equal to \eqref{Polch}. Hence it suffices to show that the elements
\eqref{SLnBasisMonomials} span the space $\widetilde{\mathbb{W}}_\br$.

Assume that there exists an element $\xi\in (\widetilde{\mathbb{W}}_\br)^*$ vanishing on all the monomials \eqref{SLnBasisMonomials}. We want to show that
in this case $\xi$ is zero; equivalently, we need to prove that $f_\xi=0$.
A non-zero polynomial divisible by $\prod_{I\prec J}(Y_{I,p_1}-Y_{J,p_2})^{k(I,J)}$ contains a monomial
\begin{equation}\label{Ymon}
\prod_{I\subset\{1,\dots,n\}} Y_I^{l_{1,I}}\dots Y_I^{l_{r_I,I}},\quad 0\le l_{1,I}\le \dots\le l_{r_I,I}
\end{equation}
such that $l_{1,I}\geq \sum_{J\prec I} r_J J$ (coming from the choice of the term $Y_{J,p_2}^{k(I,J)}$ in each bracket
$(Y_{I,p_1}-Y_{J,p_2})^{k(I,J)}$). However, the coefficient in front of monomial \eqref{Ymon} is equal to
$\xi(m_J^{(l_{1,J})}\dots m_J^{(l_{r_J,J})})$. Therefore, if $\xi$ vanishes on all the elements \eqref{SLnBasisMonomials}, then
$f_\xi$ is zero.
\end{proof}

\subsection{Arcs for toric degeneration}
We need the following Lemma.
\begin{lem}\label{SLnIdealNotIntersects}
 The polynomial algebra $\Bbbk[z_{uv}^{(k)}]_{(u,v) \neq (n,n)}$ does not intersect the defining ideal
\[\Bigg \langle det\left(
                                                 \begin{array}{ccc}
                                                   z_{11}(s) & \cdots & z_{n1}(s) \\
                                                   \vdots & \ddots & \vdots \\
                                                   z_{1n}(s) & \cdots & z_{nn}(s) \\
                                                 \end{array}
                                               \right)-1
\Bigg \rangle.\]
\end{lem}
\begin{proof}
 Rewrite the generator of ideal in the following form:
\begin{equation}
 det\left(
                                                 \begin{array}{ccc}
                                                   z_{11}(s) & \cdots & z_{n1}(s) \\
                                                   \vdots & \ddots & \vdots \\
                                                   z_{1n}(s) & \cdots & z_{nn}(s) \\
                                                 \end{array}
                                               \right)-1=
\sum_{i=1}^nz_{ni}(s)m_{\{1,\dots,n\}\backslash \{i\}}(s)-1.
\end{equation}

Each $s$-coefficient of this element is a linear function on $z_{nn}^{(k)}$ with coefficients in $\Bbbk[z_{uv}^{(k)}]_{(u,v) \neq (n,n)}$.
Moreover this is a triangular linear system, because the coefficient of $s^l$ depends only on $z_{nn}^{(l)}$, $l \leq k$.
Hence these coefficients generate a proper ideal in the polynomial algebra $\mathcal{F}[z_{nn}^{(l)}]$, where  $\mathcal{F}=\Bbbk(z_{uv}^{(k)})_{(u,v) \neq (n,n)}$.
Thus this ideal doesn't intersect $\mathcal{F}$. This completes the proof.
\end{proof}

Consider the following order  $\lhd$ on variables $z_{uv}^{(k)}$. It is lexicographic with respect to lower indices and it doesn't feel the upper indices.
More precisely $z_{uv}^{(k)}  \lhd z_{u'v'}^{(k')}$ if $v<v'$ or $v=v'$ and $u<u'$. Note that all variables are comparable
with respect to this order and for fixed $u,v$ variables $z_{uv}^{(k)}$ are equivalent.
We introduce the degree-restricted lexicographic order on monomials in $z_{uv}^{(k)}$.
This is a monomial order, i. e. for monomials $a,b,c$ if $a \lhd b$, then $ac \lhd bc$.

Note that all coefficients of minors $m_I^{(k)}$ lie in the polynomial algebra $\Bbbk[z_{uv}^{(k)}], (u,v) \neq (n,n)$ due to Lemma \ref{SLnIdealNotIntersects}.
Then the following Lemma is obvious.
\begin{lem}\label{LeadingTermMinor}
 The leading part of $m_I(s)$ with respect to $\lhd$ is equal to
\[d_I(s):=z_{i_1 1}(s)z_{i_2 2}(s) \dots z_{i_{|I|}|I|}(s).\]
\end{lem}

We denote by $d_I^{(k)}$ the coefficient of $s^k$ in the series $d_I(s)$, i. e.:
\[d_I(s):=\sum_{k=0}^\infty d_I^{(k)}.\]

\begin{rem}
The algebra $\Bbbk[d_I^{(k)}]$ is the reduced function algebra on the affine toric degeneration of
flag variety. The classical toric degeneration of the flag variety is due to Gonciulea-Lakshmibai \cite{GL}.
\end{rem}

The goal of the remainder of this Section is to prove that the leading parts of monomials from \eqref{SLnBasisMonomials}
are linear independent.

\begin{lem}\label{SLnDiagonalsLeading}
 The leading part of
\[\prod_{I\subset\{1,\dots,n\}} m_I^{(l_{1,I})}\dots m_I^{(l_{r_I,I})},\quad 0\le l_{1,I}\le \dots\le l_{r_I,I}\]
is equal to
\[\prod_{I\subset\{1,\dots,n\}} d_I^{(l_{1,I})}\dots d_I^{(l_{r_I,I})},\quad 0\le l_{1,I}\le \dots\le l_{r_I,I}\]
\end{lem}
\begin{proof}
The order $\lhd$ is monomial. Therefore the leading part of a product is equal to the product of leading parts.
\end{proof}

We have:
\[d_I^{(k)} \in \Bbbk[z_{ij}^{(k)}]_ {i \geq j}.\]
Consider the vector space
\[U_{\br}:=\langle \prod_{I\subset\{1,\dots,n\}} d_I^{(l_{1,I})}\dots d_I^{(l_{r_I,I})},\quad 0\le l_{1,I}\le \dots\le l_{r_I,I}\rangle\]
with fixed vector of numbers $\br=(r_I)$ and arbitrary $(l_{p,I})$.
We study the functional realization of the dual space $U_{\br}^*$.
We take $\sum_I r_I$ variables $Y_{I, p}$, $1 \leq p \leq r_I$ and
consider the expression
\begin{equation}\label{dDifferentParemeter}
 \prod_I d_I(Y_{I,1})\dots d_I(Y_{I,r_I}).
\end{equation}
Clearly the coefficients of this expression are all of the form $ \prod_I d_I^{(l_{1,I})}\dots d_I^{(l_{r_I,I})}$.
They span a subspace of the space spanned by elements $\prod_{u\geq v}{z_{uv}^{(k_{uv,1})}}\dots {z_{uv}^{(k_{uv,a_{uv}})}}$,
where $a_{uv}=\sum_{I|i_{v}=u} r_I$, i. e. $a_{uv}$ is equal to the number of appearance of $z_{uv}$ as a factor of $d_I$. Denote
this vector of numbers $a_{uv}$ by $\ba$ and this subspace by
$V_{\ba}$.

Consider the expression

\begin{equation}\label{zDifferentParemeter}
 \prod_{uv} z_{uv}(X_{uv,1})\dots z_{uv}(X_{uv,a_{uv}}).
\end{equation}

As before for each $\chi \in V_{\ba}^*$ we define the function:
\[f_{\chi}:=\chi\big( \prod_{uv} z_{uv}(X_{uv,1})\dots z_{uv}(X_{uv,a_{uv}})\big).\]
The functions $f_{\chi}$ form the following space (in fact, the algebra):
\[\Bbbk[X_{uv,p}]^{\mathfrak{S}_{a_{11}}\times \mathfrak{S}_{a_{21}}\times \dots}\]
of polynomials in $X_{uv,p}$ symmetric in all group of variables $X_{uv,1}, \dots, X_{uv,a_{uv}}$.
Thus the space $V_{\ba}^*$ is naturally isomorphic to this space of polynomials.

Note that we have the following substitution map
\begin{equation}\label{specializationXtoY}
\varphi_{\br}:X_{uv, \sum_{i_v=j_v=u, I \lhd J}{r_I}+j}\mapsto Y_{J,j}.
\end{equation}

Then to each variable $Y_{I,j}$ we substitute $|I|$ variables $X_{i_kk,p}$. By the construction we have:
\begin{equation}
\varphi_{\br}\left( \prod_{uv} z_{uv}(X_{uv,1})\dots z_{uv}(X_{uv,a_{uv}})\right)= \prod_I d_I(Y_{I,1})\dots d_I(Y_{I,r_I}).
\end{equation}

We have:
\[U_{\br}^*=V_{\ba}^*/ann(U_{\br}).\]
Thus the functional realization of $U_{\br}^*$ can be obtained in the following way.
We fix $\chi \in V_{\ba}*$.
Then

\[f_{\chi}(Y_{I,p})=\chi( \prod_I d_I(Y_{I,1})\dots d_I(Y_{I,r_I})).\]
Hence the dual space of $U_{\br}$ is naturally isomorphic to
$\varphi_{\br}(\Bbbk[X_{uv,p}]^{\mathfrak{S}_{a_{11}}\times \mathfrak{S}_{a_{21}}\times \dots})$.

\begin{rem}
Note that $\varphi_{\br}(\Bbbk[X_{uv,p}]^{\mathfrak{S}_{a_{11}}\times \mathfrak{S}_{a_{21}}\times \dots})\subset
\Bbbk[Y_{I,j}]^{\times_{I}\mathfrak{S}_{|I|}}$. In general this inclusion is not an equality.
\end{rem}

\begin{example}\label{SL3SpecializationExample}
Consider the case $n=3$. Then we have $6$ Pl\"ucker variables with leading monomials
\[d_1=z_{11}, d_2=z_{21}, d_3=z_{31}, d_{12}=z_{11}z_{22}, d_{13}=z_{11}z_{32}, d_{23}=z_{21}z_{32}.\]
Take a vector $\br=(r_1, r_2, r_3, r_{12}, r_{13}, r_{23})$. Then:
\[a_{11}=r_1+r_{12}+r_{13},a_{21}=r_2+r_{23}, a_{31}=r_3, a_{22}=r_{12}, a_{32}=r_{13}+r_{23}.\]

The space $V_{\ba}^*$ is isomorphic to the following algebra:

\begin{multline}
\Bbbk[X_{11,1}, \dots, X_{11, a_{11}},X_{21,1}, \dots, X_{21, a_{21}},X_{31,1}, \dots, X_{31, a_{31}},X_{22,1}, \dots, X_{22, a_{22}},\\ X_{32,1},
 \dots, X_{32, a_{32}}]^{\mathfrak{S}_{a_{11}}\times \mathfrak{S}_{a_{21}}\times\mathfrak{S}_{a_{31}}\times \mathfrak{S}_{a_{22}}\times
 \mathfrak{S}_{a_{32}}}
\end{multline}
i. e. the algebra of polynomials symmetric in $5$ groups of variables.

The map $\varphi_{\br}$ is as follows:
\[X_{11,p}\mapsto Y_{12,p};~ p=1, \dots, r_{12}, X_{11,r_{12}+p}\mapsto Y_{13,p},~ p=1, \dots, r_{13}, X_{11,r_{12}+r_{13}+p}\mapsto Y_{1,p},~ p=1, \dots, r_{1};\]
\[X_{21,p}\mapsto Y_{23,p};~ p=1, \dots, r_{23}, X_{21,r_{23}+p}\mapsto Y_{2,p};~ p=1, \dots, r_{2};\]
\[X_{31,p}\mapsto Y_{3,p};~ p=1, \dots, r_{23};\]
\[X_{22,p}\mapsto Y_{12,p};~ p=1, \dots, r_{12};\]
\[X_{32,p}\mapsto Y_{13,p};~ p=1, \dots, r_{13}, X_{32,r_{13}+p}\mapsto Y_{23,p};~ p=1, \dots, r_{23}.\]

In particular the image of $\varphi_{\br}$ is the subalgebra in $\Bbbk[Y_{I,p}]$ generated by $5$ subalgebras isomorphic to symmetric algebras
in $a_{ij}$ variables.
\end{example}

The next Lemmas contain the properties of symmetric and supersymmetric polynomials which we need for the proof.

\begin{dfn}\label{PartialSymmetricAlgebra}
We take $k$ groups of variables $Y_{l,i}$, $l=1, \dots, k$, $i=1, \dots, r_l$ and $q$ subsets $S_j \subset\{1, \dots, k\}$.
We denote:
\[\mathcal{R}_{S_j}=\Bbbk[Y_{l,i}]_{l \in S_j, ~i=1, \dots, r_l}^{\mathfrak{S}_{\sum_{l \in S_j}r_{l}}}.\]
We define the subalgebra generated by $k$ subalgebras of symmetric functions:
\[\mathcal{R}_{\{S_1, \dots, S_q\}}:=\Bbbk[\mathcal{R}_{S_1}, \dots, \mathcal{R}_{S_k}]\subset \Bbbk[Y_{l,i}]_{l=1, \dots, k, ~i=1, \dots, r_l}.\]
\end{dfn}

\begin{lem}\label{SymmetricGroups}
Assume that $q=k$ and the indicator matrix
$(\delta_i^j)$, $\delta_i^j=1$ if $i \in S_j$, $0$ otherwise, is invertible. Then:
\[\mathcal{R}_{\{S_1, \dots, S_q\}}=\Bbbk[Y_{l,i}]^{\mathfrak{S}_{r_1}\times\dots \times \mathfrak{S}_{r_k}}_{l=1, \dots, k, ~i=1, \dots, r_l}.\]
\end{lem}
\begin{proof}
It is sufficient to prove that every power sum of elements in each group belongs to this algebra, i. e.:
\[p_{l,m}=\sum_{i=1}^{r_l} Y_{l,i}^m\in \Bbbk[\mathcal{R}_{S_1}, \dots, \mathcal{R}_{S_k}].\]
However by the construction we have:
\[\sum_{l \in S_j} p_{l,m}\in \mathcal{R}_{S_j} \subset \Bbbk[\mathcal{R}_{S_1}, \dots, \mathcal{R}_{S_k}].\]

Then using the invertibility of $(\delta_i^j)$ we have that $p_{l,m}$ are linear combinations of these elements.
This completes the proof of Lemma.
\end{proof}

The next Lemma is a property of supersymmetric functions. Recall the definition.

Take two groups of variables $A_1, \dots, A_r$, $B_1, \dots, B_l$. A polynomial $f \in \Bbbk[A_i, B_j]^{\mathfrak{S}_r\times \mathfrak{S}_l}$
is called supersymmetric if the following holds for the substitution $A_r \mapsto X, B_l \mapsto X$:
\[f|_{A_r \mapsto X, B_l \mapsto X}\in \Bbbk[A_i, B_j, i\leq r-1, j \leq l-1]^{\mathfrak{S}_{r-1}\times \mathfrak{S}_{l-1}}.\]
I. e. this substitution doesn't depend on the variable $X$.
\begin{lem}\label{superSymResultant}
The  ring of supersymmetric functions is generated by coefficients of the series
\[\frac{\prod_{j=1}^r(1-tA_j)}{\prod_{j=1}^l(1-tB_j)}.\]
The ideal of the ring $\Bbbk[A_1, \dots, A_r, B_1, \dots, B_l]^{\mathfrak{S}_r\times \mathfrak{S}_l}$ generated by
\begin{equation}\label{ResultantAB}
  \prod_{1 \leq i \leq r, 1 \leq j \leq l} (A_i-B_j)
\end{equation}
belongs to the ring of supersymmetric polynomials and is equal to the kernel of the substitution map $f \mapsto f|_{A_r \mapsto X, B_l \mapsto X}$.
\end{lem}
\begin{proof}
The first claim is proved in \cite{M}, Chapter 1.5. The second claim is obvious.
\end{proof}

Take two numbers $a_{\alpha}$, $a_{\beta}$
and two groups of variables $X_{\alpha,1}, \dots, X_{\alpha,a_{\alpha}}$, $X_{\beta,1}, \dots, X_{\beta,a_{\beta}}$.
Define the vector of numbers $\br=(r_b, r_c, r_d)$ such that $r_b+r_c=a_{\alpha}$, $r_d+r_c=a_{\beta}$
and take three groups of variables $B_{1}, \dots, B_{{r_b}}$, $C_{1}, \dots, C_{{r_c}}$, $D_{1}, \dots, D_{{r_d}}$.
Consider the map $\varphi_\br:\Bbbk[X_{\alpha,1}, \dots, X_{\alpha,a_{\alpha}}$, $X_{\beta,1}, \dots, X_{\beta,a_{\beta}}]^{\mathfrak{S}_{a_{\alpha}}
\times \mathfrak{S}_{a_{\beta}}} \rightarrow \Bbbk[B_{1}, \dots, B_{{r_b}}, C_{1}, \dots, C_{{r_c}}, D_{1}, \dots, D_{{r_d}}]^{
\mathfrak{S}_{r_b}\times \mathfrak{S}_{r_c}\times \mathfrak{S}_{r_d}}$:
\[\varphi_{\br}(X_{\alpha,1})=B_{1},\dots, \varphi_{\br}(X_{\alpha,a_b})=B_{a_b},
\varphi_{\br}(X_{\alpha,a_b+1})=C_{1},\dots, \varphi_{\br}(X_{\alpha,a_b+a_c})=C_{a_c};\]
\[\varphi_{\br}(X_{\beta,1})=D_{1},\dots, \varphi_{\br}(X_{\beta,a_d})=D_{a_d},
\varphi_{\br}(X_{\beta,a_d+1})=C_{1},\dots, \varphi_{\br}(X_{\beta\alpha,a_d+a_c})=C_{a_c}.\]
We consider the order $\br<\br'=(r'_b, r'_c, r'_d)$ if $r_c<r'_c$.
\begin{lem}\label{subspaceMultResultant}
There exists a subspace $S_\br \subset \Bbbk[X_{\alpha,1}, \dots, X_{\alpha,a_{\alpha}}$, $X_{\beta,1}, \dots, X_{\beta,a_{\beta}}]^{\mathfrak{S}_{a_{\alpha}}
\times \mathfrak{S}_{a_{\beta}}} $
such that $\varphi_{\br'}(S_\br)=\{0\}$ for $\br'>\br$ and
\[\varphi_{\br}(S_\br)=\prod_{1 \leq i \leq r_b, 1 \leq j \leq r_d}(B_i-D_j)\Bbbk[B_{1}, \dots, B_{{r_b}}, C_{1}, \dots, C_{{r_c}}, D_{1}, \dots, D_{{r_d}}]^{
\mathfrak{S}_{r_b}\times \mathfrak{S}_{r_c}\times \mathfrak{S}_{r_d}}.\]
\end{lem}
\begin{proof}
We denote
\[\mathcal{R}(\br):=\Bbbk[B_{1}, \dots, B_{{r_b}}, C_{1}, \dots, C_{{r_c}}, D_{1}, \dots, D_{{r_d}}] ^{
\mathfrak{S}_{r_b}\times \mathfrak{S}_{r_c}\times \mathfrak{S}_{r_d}}.\]
Define a substitution map $\mathcal{R}(\br)\rightarrow \mathcal{R}(\br')$ for $\br'>\br$.
\[\eta_{\br}^{\br'}: B_{r_b}\mapsto C_{r_c+1},D_{r_d}\mapsto C_{r_c+1},\dots, B_{r_b-r'_c+r_c+1}\mapsto C_{r'_c},D_{r_d-r'_c+r_c+1}\mapsto C_{r'_c}.\]
Note that $\varphi_{\br'}=\eta_{\br}^{\br'}\circ\varphi_\br$.
Consider the coefficients of the formal series
\[\frac{\prod_{i=1}^{a_\alpha} (1-tX_{\alpha,i})}{\prod_{i=1}^{a_\beta} (1-tX_{\beta,i})}.\]
By Lemma \ref{superSymResultant} their image under $\varphi_\br$ generate the ring of supersymmetric functions in
$B_{1}, \dots, B_{{r_b}}$ and $D_{1}, \dots, D_{{r_d}}$ and in particular the space
\[\prod_{1 \leq i \leq r_b, 1 \leq j \leq r_d}(B_i-D_j)\Bbbk[B_{1}, \dots, B_{{r_b}}, D_{1}, \dots, D_{{r_d}}]^{
\mathfrak{S}_{r_b}\times \mathfrak{S}_{r_d}}.\]
Using Lemma \ref{PartialSymmetricAlgebra} we have that the image of $\varphi_\br$ contains the whole principal ideal generated by
$\prod_{1 \leq i \leq r_b, 1 \leq j \leq r_d}(B_i-D_j)$.
However:
\[\eta_{\br}^{\br'}\prod_{1 \leq i \leq r_b, 1 \leq j \leq r_d}(B_i-D_j)=0.\]
Therefore the preimage of this subspace under $\varphi_\br$ is the needed subspace $S_\br$. This completes the proof of Lemma.
\end{proof}

Fix a vector $\ba=(a_{uv})_{u \geq v}$.
We recall the order $\prec$ on the vectors $\br=(r_I)$.
Our goal is to study the symmetric polynomials in $X_{uv,i}$ which are annihilated by all maps $\varphi_{\br'}$ for
$\br'\succ\br$.

We use the following notation to simplify the formulas:
\begin{equation}\label{SimpleProduct}
\prod_{1 \leq i \leq r_{I},1 \leq j \leq r_{J}}(Y_{I,i}-Y_{J,j})=:\big(Y_{I,i}-Y_{J,j}\big).
\end{equation}
We next prove the needed property of symmetric functions in the case $r_I=0$ for $|I|\geq 3$ or equivalently $a_{uv}=0$ for $v \geq 3$.
We omit the braces when we write the set $I$, i. e. we use a notation $\{a,b\}=ab$. 
\begin{prop}\label{phiImageTworows}
Assume that $a_{uv}=0$ for $v \geq 3$. Then for each $\br$ there exists a subspace $S_\br$ such that
\[\varphi_{\br'}(S_\br)=0 ~\text{for}~\br'\succ\br;\]
\[\varphi_{\br}(S_\br)=
\prod_{a<b<c<d}(Y_{ad,i}-Y_{bc,j})\Bbbk[Y_{I,j}]^{\times_{I}\mathfrak{S}_{|I|}}.
\]
\end{prop}
\begin{proof}
For a subset $I$ we define the set of variables $X_{uv,l}^I$, $l=1, \dots, a_{uv}^I$, where
\[a_{uv}^I=a_{uv}-\sum_{j_v=u,J\preceq I}r_J.\]
Define then the algebra $\mathcal{R}^I$ of polynomials symmetric in groups of variables $Y_{J,l}$, $J \prec I$ and
$X_{uv,l}^I$. Note that the order $\prec$ is linear, therefore the largest element $I'$ such that
$I' \prec I$ is well defined, denote it by $I^{-}$.
In the same way we denote by $I^+$ the smallest element larger then $I$.
We denote $I^{++}:=(I^+)^+$, etc.
 Then we have the map $\varphi_{\br}^I:\mathcal{R}^{I^-}\rightarrow \mathcal{R}^I$:
 \[\varphi_{\br}^I(X_{l}^{I^-})=
Y_{I,l},~ l \leq r_I;~\varphi_{\br}^I(X_{uv,r_I+l}^{I^-})=
X_{uv,l}^{I}, ~\text{if}~ I_v=u;\]
\[\varphi_{\br}^I(X_{uv,l}^{I^-})=X_{uv,l}^{I}, ~\text{otherwise};\]
\[\varphi_{\br}^I(Y_{J,l})=Y_{J,l}, J \prec I.\]
Then $\varphi_{\br}$ is equal to the composition of all maps $\varphi_{\br}^I$.

Note that the conditions $|I|\leq 2$ means that if $(u,v) \not \in \{(i_1,1),(i_2,2)\}$, then the map
 $\varphi_{\br}^I$ bijectively sets the group of variables $X_{uv,p}^{I^-}$ to the group of
 variables $X_{uv,p}^{I}$.
Thus we can apply Lemma \ref{subspaceMultResultant} to the groups of variables $X_{i_11,p}^{I^-}$,
$X_{i_22,p}^{I^-}$.
Hence there exists the
subspace $S_\br^I$ whose image is the principal ideal generated by
\[\prod_{l=1,\dots,a_{i_11}^I,m=1,\dots,a_{i_22}^I}(X_{i_11,l}^{I}-X_{i_22,m}^{I})\]
and $\varphi_{\br'}^I(S_\br^I)=0$ if $r'_I>r_I$.

One can compute:
\[\varphi_{\br}^{\{n\}}\circ \varphi_{\br}^{\{n-1\}} \circ\dots \circ \varphi_{\br}^{I^{++}}\circ\varphi_{\br}^{I^+}
\left(\prod_{l=1,\dots,a_{i_11}^I,m=1,\dots,a_{i_22}^I}(X_{i_11,l}^{I}-X_{i_22,m}^{I})\right)=
\prod_{i_1<b<i_2<c}\big(Y_{i_1c}-Y_{bi_2}\big).\]

Thus comparing the properties of maps $\varphi_{\br}^I$ we obtain that there is a subspace $S_\br$ such that for each $I$
its image under $\varphi_\br^{I^{-}}\circ \dots\varphi_\br^{12}$ is contained in $S_\br^I$ and $\varphi(S_\br)$ is equal to
$\prod_{a<b<c<d}(Y_{ad,i}-Y_{bc,j})\Bbbk[Y_{I,j}]^{\times_{I}\mathfrak{S}_{|I|}}$. If $\br'\succ\br$, then
there exists the smallest $I$ such that $r_I<r'_I$. Then the algebras $\mathcal{R}^J$ are the same for $J\prec I$.
Therefore the space $S_\br^I$ is contained in the kernel of $\varphi_{\br'}^I$ and thus $S_\br$ is contained in
the kernel of $\varphi_{\br'}$. This completes the proof of Proposition.
\end{proof}

Now we a ready to prove the needed property of the morphism $\varphi_\br$ in the whole generality.

\begin{thm}\label{phiImage}
For each $\br$ there exists a subspace $S_\br\subset \Bbbk[X_{uv,i}]^{\times \mathfrak{S}_{a_{uv}}}$ such that
\[\varphi_{\br'}(S_\br)=0 ~\text{for}~\br'\succ\br;\]
\[\varphi_{\br}(S_\br)=
\prod_{I\prec J}(Y_{I,i}-Y_{J,j})^{k(I,J)}\Bbbk[Y_{I,j}]^{\times_{I}\mathfrak{S}_{|I|}}.
\]
\end{thm}
\begin{proof}
We fix a number $l$, $1 \leq l \leq n-1$. Recall the truncation map. For $I \subset \{1, \dots, n\}$, $|I| \geq l$, we put
\[\tr^{l-1}(I):=\{i_{l}, \dots,i_{|I|}\},\]
i. e. this operation deletes from $I$ its $l-1$ smallest elements and is $l-1$st power of the operator $\tr$ deleting the smallest element of the set.

For each subset $U \subset \{l,\dots,n-1\}$ we define
\[r_U^{l-1}:=\sum_{\tr^{l-1}(I)=U}r_I.\]
Denote the vector of these numbers by $\br^{l-1}$.

We construct the algebra $\mathcal{R}_{\br^{l-1}}^{l-1}$ as a symmetric algebra in $X_{uv,p}, v \leq l-1, p=1, \dots a_{uv}$ and
$Y_{U,p}^{l-1}$, $U \subset \{l,\dots,n-1\}$, $1 \leq p \leq r_U^{l-1}$. In particular $\mathcal{R}_{\br}^{0}$
is a symmetric algebra is $Y_{I,p}$.

Define a map $\psi_\br^l:\mathcal{R}_{\br^{l}}^{l}\rightarrow \mathcal{R}_{\br^{l-1}}^{l-1}$:
\[Y_{U,p}^{l}\mapsto Y_{\{1\} \cup U,p}^{l-1}, p=1, \dots, r_{\{1\} \cup U}^{l-1};
Y_{U,r_{\{1\} \cup U}^{l-1}+p}^{l}\mapsto Y_{\{1\} \cup U,p}^{l-1}, p=1, \dots, r_{\{2\} \cup U}^{l-1}, \dots; \]
\[X_{ul,\sum_{U' \prec U}r_{U'}^{l} p} \mapsto Y_{\{u\}\cup U,p}.\]
Then we have
\[\varphi_\br:=\psi_\br^1\circ \psi_\br^2 \dots \circ \psi_\br^{n-1}.\]
We equip the elements $Y_{U,p}^{l}$ the lexicographic order, i. e. $U \blacktriangleleft V$ if
for some $p$ their $p-1$ smallest elements coincide and the $p$th element of $U$ is smaller then the $p$th element on $V$. Note
that this order doesn't coincide with the order $\prec$.
Then we define the lexicographic order on the vectors $\br^{l-1}$ with the same truncation to $\br^{l}$.
Thus using Proposition \ref{phiImageTworows} we have that there exists a subspace $S_{\br^{l-1}}\subset \mathcal{R}_{\br^{l}}^{l}$
such that for each $\br'^{l-1}\succ\br^{l-1}$ $\psi_{\br'^{l-1}}^{l-1}(S_{\br^{l}})=0$ and
\[\psi_{\br'^{l-1}}^{l-1}(S_{\br^{l-1}})=\prod_{U\blacktriangleleft V, 1 \leq a<b\leq n}\big( Y_{U\cup \{b\}}-Y_{V\cup \{a\}} \big)
\mathcal{R}_{\br^{l-1}}^{l-1}.\]
Then we complete the proof of Theorem by induction.
\end{proof}

\begin{cor}\label{linearIndependentMonomials}
The monomials \eqref{SLnBasisMonomials} as elements of $\mathbb{W}$ are linear independent.
\end{cor}
\begin{proof}
The linear span of leading parts of monomials \eqref{SLnBasisMonomials} has the same graded dimension as the
degenerate algebra $\widetilde{\mathbb{W}}$. However the graded dimension of the degenerate algebra is greater then or equal
to the graded dimension of original algebra and the linear span of leading parts of a set of polynomials is less then or
equal to the graded dimension of the linear span of this set of polynomials. Therefore these dimensions are equal and the set
of monomials \eqref{SLnBasisMonomials} is linear independent. This completes the proof of Corollary.
\end{proof}

\section{Semi-infinite Pl\"ucker relations in type $C$}
We consider the current group of the symplectic group. The algebra of functions $\Bbbk[SP_{2n}[[t]]]$ on this group is generated by elements $z_{uv}^{(k)}$,
$u,v \in \{1,\bar 1, \dots, n, \bar n\}$, $k=0,1,\dots$ and is defined by the relations:
\[\sum_{l=1}^n z_{lu}(s)z_{\bar{l}v}(s)-z_{\bar{l}u}(s)z_{lv}(s)=
\begin{cases}
  1, & \mbox{if } v=\bar{u} \\
  -1, & \mbox{if } u=\bar{v} \\
  0, & \mbox{otherwise}.
\end{cases}\]

Analogously to the finite-dimensional case we have:
\begin{prop}\label{SPCurrentlinearDependenceMinors}
For any $I \subset \{1, 2, \dots, n, \bar {1}, \dots, \bar{n}\}$, $|I| \leq n-2$
we have:
\[m_{I \cup \{1, \bar{1}\}}(s)+m_{I \cup \{2, \bar{2}\}}(s)+\dots++m_{I \cup \{n, \bar{n}\}}(s)=0,\]
where union means the union of multisets.
\end{prop}

Our plan is to study the subalgebra of $\Bbbk[SP_{2n}[[t]]]$ invariants under the right $\fn_+[t]$ action.
We first need to prove the analogue of Proposition \ref{SLnCurrentInvariants} for the symplectic algebra.

\begin{prop}\label{SPnCurrentInvariants}
\[\Bbbk[SP_n[[t]]]^{\fn_+[t]}=
\Bbbk[m_I^{(k)}]_{I ~\text{is allowed}}\]
\end{prop}
\begin{proof}
The proof literally coincides with the proof of Proposition \ref{SLnCurrentInvariants}.


\end{proof}

It is clear that the minors in type $C_{n}$ satisfy the same (semi-infinite Pl\"ucker) relations as the minors in type
$A_{2n}$. More precisely for allowed $I,J$ we define the number $k(I,J)$ as in Definition \ref{kdefinition}.

\begin{prop}\label{SPpropsnakeplueckerequation}
$a)$.\ For any $k' \leq k(I,J)-1$ we have the following equality in $\mathbb{W}[[s]]$:
\begin{equation}\label{SPsnakeplueckerequation}
\sum_{A\subset S(I,J), |A|=|I \cap S|}(-1)^{{\rm sign}(A)}
\frac{\partial^{k'} m_{I \backslash S \cup A}(s)}{\partial s^{k'}} m_{I \backslash S \cup (S \backslash A)}(s)=0.
\end{equation}
$b).$\ We have
\[\frac{\partial^{k'} m_{I}(s)}{\partial s^{k'}}m_{J}(s)\]
is the leading part of the previous equation.
\end{prop}

Thus one can define the degenerated algebra $\widetilde{\mathbb{W}}$ generated by allowed minors $m_I^{(k)}$ and satisfying the relations
\begin{equation}\label{SPdegeneratedEquation}
\frac{\partial^{k'} m_{I}(s)}{\partial s^{k'}}m_{J}, ~k'<k(I,J).
\end{equation}
The graded character of this algebra is greater then or equal to the graded character of $\mathbb{W}$.
Attach the degree $I+l \delta$ to the allowed monomial $m_I^{(l)}$.
\begin{prop}\label{SPdegeneratedInclusion}
The algebra $\widetilde {\mathbb{W}}$ is graded by the free group in generators $I \subset \{1,\bar 1 \dots, n, \bar n\}$, $i_p \geq p$ and
symbol $\delta$.
It is isomorphic to the subalgebra of the degenerated algebra $\widetilde {\mathbb{W}}$ of type $A_{2n}$
generated by minors $m_I^{(l)}$, $i_p \geq 2p-1$.
In particular the character of the homogenous subset
$\widetilde {\mathbb{W}}_\br \subset \widetilde {\mathbb{W}}$, $\br= \sum r_I I$, of type $C_n$ is equal to
\begin{equation}\label{SPdegeneraterBasis}
\ch \widetilde{\mathbb{W}}_\br=
\frac{q^{\sum_{I \prec J}k(I,J)r_Ir_J}}
{\prod_{I} (q)_{r_I}}.
\end{equation}
The homogeneous component $\widetilde{\mathbb{W}}_\br \subset \widetilde{\mathbb{W}}$ has the basis consisting of monomials of the form
\begin{equation}\label{SPnBasisMonomials}
\prod_{I\subset\{1,\dots,n\}} m_I^{(l_{1,I})}\dots m_I^{(l_{r_I,I})},\quad 0\le l_{1,I}\le \dots\le l_{r_I,I}
\end{equation}
such that $l_{1,I}\geq \sum_{J \prec I} k(I,J) r_J$, $m_I, m_J$ are allowed minors.
\end{prop}
\begin{proof}
The defining relations are homogenous with respect to the grading coming from the attaching degree $I+l\delta$ to $m_I^{(l)}$.

We identify the algebra $\widetilde{\mathbb{W}}$ of type $C_n$ with subalgebra of the algebra $\widetilde{\mathbb{W}}$ of type $A_{2n}$
in the following way.
We identify the subsets $I \subset \{1,\bar 1 \dots, n, \bar n\}$, $i_p \geq p$
with the subsets $I \subset \{1, \dots, 2n\}$, $i_p \geq 2p-1$ by the map on elements $p \mapsto 2p-1$, $\bar p \mapsto 2p$.
The map is well defined because the defining relations of the algebra $\widetilde{\mathbb{W}}$ of type $C_n$ coincide
with the defining relations of this subalgebra.

Each graded component $\widetilde{\mathbb{W}}_\br$ goes to the graded component and
the basis \ref{SPnBasisMonomials} of this component goes to the basis \ref{SLnBasisMonomials}
In particular the image of this inclusion is the sum of some homogenous components of the graded algebra.
This gives the needed character formula.
\end{proof}

Next we prove the current analogue of Lemma \ref{idealRelationsIntersection}
\begin{lem}\label{SPcurrentidealRelationsIntersection}
The defining ideal of $\Bbbk[SP_{2n}[[t]]]$ has trivial intersection with the polynomial subalgebra $\Bbbk[z_{uv}^{(k)}]$, $v \in \{1, 2, \dots, n\}$,
$(u,v) \neq (\bar p, l), p < l$.
\end{lem}
\begin{proof}
The proof is parallel to the proof of Lemma \ref{idealRelationsIntersection}. Consider first the algebra
generated by elements $z_{uv}^{(k)}$, $v \in \{1, 2, \dots, n\}$;
\[\mathcal{G}':=\Bbbk[z_{uv}^{(k)}]/\langle \sum_{k=1}^n z_{uv}(s)z_{\bar{u}v}(s)-z_{\bar{u}v}(s)z_{uv}(s)\rangle.\]
Take the field $\mathcal{F}=\Bbbk(z_{uv}^{(0)})$, $v \in \{1, 2, \dots, n\}$,
$(u,v) \neq (\bar p, l), p < l$.
Define the algebra
\[\mathcal{G}:=\mathcal{F}[z_{uv}^{(k)}]/\langle \sum_{k=1}^n z_{uv}(s)z_{\bar{u}v}(s)-z_{\bar{u}v}(s)z_{uv}(s)\rangle,\]
$(u,v,k) \neq ((\bar p, l, 0), p < l$.
This is the algebra of fractions for $\mathcal{G}'$.
The series $z_{uv}(s)\subset \mathcal{G}[[s]]$, $(u,v) \neq (\bar p, l), p < l$,  is invertible (because its free term is invertible).
Then  as in the proof of Lemma \ref{idealRelationsIntersection} we have that the set of $n(n-1)/2$ equations
\[\sum_{v=1}^n z_{uv}(s)z_{\bar{u}v}(s)-z_{\bar{u}v}(s)z_{uv}(s)=0\]
is a triangular linear system for the $n(n-1)/2$ variables $z_{\bar p, l}(s), p < l$ and the diagonal coefficients of this system are invertible.
Therefore the defining relations of $\mathcal{G}$ can be rewritten in the form:
\[z_{\bar p l}(s)-a_{\bar p l}(s)=0, p < l, a_{\bar p l}\in \mathcal{F}[z_{uv}^{(k)}], v \in \{1, 2, \dots, n\}, (u,v) \neq (\bar p, l), p < l, k \geq 1.\]
These are the linear equations one for each element $z_{\bar p l}^{(k)}, p < l$. Hence the ideal generated by these equations
doesn't intersect the ring $\mathcal{F}[z_{uv}^{(k)}], v \in \{1, 2, \dots, n\}, (u,v) \neq (\bar p, l), p < l, k \geq 1$ and moreover the ring
of polynomials in $z_{uv}^{(k)}$, $v \in \{1, 2, \dots, n\}$,
$(u,v) \neq (\bar p, l), p < l$.
In particular we have $\mathcal{G} \simeq \mathcal{F}[z_{uv}^{(k)}]$, $v \in \{1, 2, \dots, n\}$,
$(u,v) \neq (\bar p, l), p < l$, $k \geq 1$.
Let $\widetilde{\mathcal{G}}$ be the fraction field of this ring. Consider the algebra over this field
generated by variables $z_{uv}^{(k)}$, $v \in \{\bar 1, \bar 2, \dots, \bar n\}$ modulo the remaining defining relations of the
algebra of functions over symplectic group, i. e. the relations that the vectors $(z_{u\bar v}(s)), u=1, \dots, n, \bar 1, \dots ,\bar n$
are dual to the vectors $(z_{u v}(s)), u=1, \dots, n, \bar 1, \dots ,\bar n$. As in the proof of Lemma \ref{idealRelationsIntersection}
the ideal generated by these relations doesn't intersect the field $\widetilde{\mathcal{G}}((s))$ because each Lagrangian subspace in
$2n$ dimensional symplectic space has a dual subspace. This completes the proof.
\end{proof}
Next define the order of elements $z_{uv}^{(k)}$, $u\in \{1, \dots, n, \bar 1, \dots, \bar n\}$, $v \in \{1, \dots, n\}$, $k \geq 0$ as the
lexicographic order in $(u,v)$, i. e. $z_{uv}^{(k)}\leq z_{u'v'}^{(k')}$, if $v < v'$ or $v=v'$ and $u\leq u'$
(recall the order $1 < \bar 1 < \dots < n < \bar n$).
Define the monomial order on the monomials in these elements as the degree restricted lexicographic order.
 Then we easily have that the leading part of the allowed minors in $m_I(s)$ is equal to
 \[d_I(s)=z_{i_1,1}(s)z_{i_2,2}(s)\dots, z_{i_{|I|},|I|}(s).\]
 Our goal is to study the algebra
 \[\Bbbk[d_I^{(k)}]\]
 of leading terms of the algebra generated by allowed minors.

\begin{prop}
The leading parts of the monomials \eqref{SPnBasisMonomials} are linear independent. In particular the monomials
\eqref{SPnBasisMonomials} are linear independent.
\end{prop}
\begin{proof}
Consider the identification $z_{uv}\mapsto z_{2u-1,v}$, $z_{\bar u v}\mapsto z_{2u,v}$. This map sends the leading terms
of monomials in type $C$ to leading parts of monomials in type $A$. Then the comparison to the similar property of
leading monomials in type $A$ completes the proof.
\end{proof}
We summarize the results of this Section in the following theorem. We use the notation $e^{\varepsilon_{\bar l}}=e^{-\varepsilon_l}$
\begin{thm}
The algebra $\mathbb{W}$ in type $SP_{2n}$ is generated by allowed minors $m_I^{(k)}$ with
defining relations \eqref{SPsnakeplueckerequation} and its basis is \eqref{SPnBasisMonomials}.
In particular the character or the Weyl module $\mathbb{W}(\sum \lambda_p \omega_p)$ is the following:
\[\sum_{\br} \frac{q^{\sum_{I \prec J}k(I,J)r_Ir_J}}
{\prod_{I} (q)_{r_I}}(\prod_{I}\prod_{k=1}^{|I|} e^{\varepsilon_{i_k}})^{r_I},\]
where the summation is on $\br$ such that $\sum_{|I|=p} r_I=\lambda_p$.
\end{thm}
\begin{proof}
We sum up the $q$ degrees of homogenous components $\widetilde{\mathbb{W}}$ multiplied by the weight of this components
with respect to Cartan subalgebra.
\end{proof}
\section*{Acknowladges}
The author is grateful to I. Makhlin and A. Khoroshkin for useful discussions, to E. Feigin for useful discussions and
the help in preparing the text.


\begin{thebibliography}{99}
\bibitem[B]{B}
A. Borel, {\it Linear algebraic groups}, Springer, 1969.


\bibitem[BC]{BC}
M.~Brion, S.~Kumar, {\it Frobenius Splitting Methods in Geometry and Representation Theory,}
Springer, 2004, PM, volume 231.

\bibitem[BF1]{BF1}
A.Braverman, M.Finkelberg, {\it Weyl modules and $q$-Whittaker functions}, Math. Ann., vol. 359 (1),
2014, pp. 45--59.

\bibitem[BF2]{BF2}
A.Braverman, M.Finkelberg, {\it Twisted zastava and $q$-Whittaker functions},
J. Lond. Math. Soc. (2) 96 (2017), no. 2, 309--325.

\bibitem[BF3]{BF3} A.Braverman, M.Finkelberg, {\it Semi-infinite Schubert varieties and quantum K-theory of flag manifolds},
J. Amer. Math. Soc. 27(2014), no. 4, 1147--1168.



\bibitem[CFK]{CFK}
V.~Chari, G.~Fourier, and T.~Khandai.
A categorical approach to Weyl modules.
Transform. Groups, 15(3):517--549, 2010.


\bibitem [CL1]{CL1}
{V.~Chari}, {S.~Loktev},
{\it Weyl, Demazure and fusion modules for
the current algebra of $\msl_{r+1}$}, Adv. Math. 207 (2006), 928--960.

\bibitem [CL2]{CL2}
{V.~Chari}, {S.~Loktev},
{\it An application of global Weyl modules of $\msl_{n+1}[t]$ to invariant theory},
Journal of Algebra, Volume 349, Issue 1, 2012, pp. 317--328.



\bibitem [CP]{CP}
{V.~Chari}, {A.~Pressley},
{\it Weyl Modules for Classical and Quantum Affine Algebras}, Represent. Theory 5 (2001),
191–223.

\bibitem [DC]{DC}
C. De Concini, {\it Symplectic standard tableaux}, Adv. Math., 34, 1--27, 1979.


\bibitem[F]{F}
W.~Fulton, {\it Young Tableaux, with Applications to Representation Theory and Geometry.}
Cambridge University Press, 1997.


\bibitem[Fr]{FBZ}
E. Frenkel, {\it Langlands correspondence for loop groups},
Cambridge Studies in Advanced Mathematics, 103. Cambridge University Press, Cambridge, 2007.

\bibitem[FiMi]{FiMi}
M.Finkelberg, I.Mirkovii\'c, {\it Semi-infinite flags I. Case of global curve $\bP^1$}.
In Differential topology, infinite-dimensional Lie algebras, and applications, volume 194 of
Amer. Math. Soc. Transl. Ser. 2, pages 81--112. Amer. Math. Soc., Providence, RI, 1999.

\bibitem[FKM]{FKM}
E. Feigin, A. Khoroshkin, I. Makedonskyi,
{\it Peter-Weyl, Howe and Schur-Weyl theorems for current groups}, arXiv:1906.03290.

\bibitem[FeMa2]{FeMa2}
E. Feigin, I. Makedonskyi, {\it Semi-infinite Plücker relations and Weyl modules},
International Mathematics Research Notices, rny121, 2019.

\bibitem[FMO]{FMO}
E.~Feigin, I.~Makedonskyi, D.~Orr, {\it Generalized Weyl modules and nonsymmetric q-Whittaker functions},
Adv. Math. 330 (2018), 997--1033.

\bibitem[GL]{GL}
N.~Gonciulea, V.~Lakshmibai {\it Degenerations of flag and Schubert varieties to toric varieties}
Transformation Groups vol. 1 (1996), 215–248.


\bibitem[I]{I} B.~Ion, {\it Nonsymmetric  Macdonald  polynomials  and  Demazure  characters},
Duke Math. J.116(2003), no. 2, 299--318.


\bibitem[Kat1]{Kat}
S.~Kato, {\it Demazure character formula for semi-infinite flag manifolds}, Math. Ann., 2018, 1-33.


\bibitem[Kn]{Kn}
A.~Knapp, {\em Lie Groups Beyond an Introduction}, 2002, Progress in Mathematics, 140 (2nd ed.), Boston: Birkh\"auser.


\bibitem[Kum1]{Kum1}
S.~Kumar, {\it Kac-Moody groups, their flag varieties and representation theory},
Progr. Math., 204. Birkh{\" a}user Boston, Inc., Boston, MA, 2002.


\bibitem[Kum2]{Kum2}
S.~Kumar, {\it Conformal Blocks, Generalized Theta Functions and Verlinde Formula}, Chapter 1, in preparation.


\bibitem[M]{M}
I.~Macdonald, {\em Symmetric functions and Hall polynomials}, 2nd ed. 1995.

\bibitem[MS]{MS}
E. Miller, B. Sturmfels {\it Combinatorial Commutative Algebra}, Springer, 2005.

\bibitem[N]{N}
K.Naoi, {\it Weyl modules, Demazure modules and finite crystals for non-simply laced type},
Adv. Math. 229 (2012), no. 2, 875--934.

\bibitem[Na]{Na} J.F.Nash,  {\it Arc structure of singularities}, Duke Math. J., 81 (1995), pp. 31–38

\bibitem[Oo]{Oo}
F.~Oort, {\it Algebraic group schemes in characteristic zero are reduced}, Invent. Math. 2 (1966), 79–80.

\bibitem[S]{S}
Y.~Sanderson, {\em On the Connection Between Macdonald Polynomials and Demazure Characters},
J. of Algebraic Combinatorics, {11} (2000), 269--275.



\end{thebibliography}
\end{document}